\documentclass[11pt]{article}
\usepackage{graphicx}
\usepackage{color}
\usepackage{amsmath}
\usepackage{amssymb}
\usepackage{amscd}
\usepackage{bbm}
\newcommand{\R}{\mathbb{R}}
\newcommand{\inr}[1]{\bigl< #1 \bigr>}

\newcommand{\E}{\mathbb{E}}

\newcommand{\eps}{\varepsilon}

\newtheorem{Theorem}{Theorem}[section]
\newtheorem{Lemma}[Theorem]{Lemma}

\newtheorem{Definition}[Theorem]{Definition}

\newtheorem{Corollary}[Theorem]{Corollary}
\newtheorem{Remark}[Theorem]{Remark}

\newtheorem{Assumption}{Assumption}[section]

\numberwithin{equation}{section}

\def \proof {\noindent {\bf Proof.}\ \ }

\def \endproof
{{\mbox{}\nolinebreak\hfill\rule{2mm}{2mm}\par\medbreak}}

\begin{document}
\title{Upper bounds on product and multiplier empirical processes}
\author{Shahar Mendelson\footnote{Department of
Mathematics, Technion, I.I.T, Haifa 32000, Israel
\newline
{\sf email: shahar@tx.technion.ac.il}
\newline
Supported in part by the Mathematical Sciences Institute,
The Australian National University, Canberra, ACT 2601,
Australia. Additional support was given by the Israel Science Foundation.
 } }

\medskip
\maketitle
\begin{abstract}
We study two empirical process of special structure: firstly, the centred multiplier process indexed by a class $F$,  $f \to \left|\sum_{i=1}^N (\xi_i f(X_i) - \E \xi f)\right|$, where the i.i.d. multipliers $(\xi_i)_{i=1}^N$ need not be independent of $(X_i)_{i=1}^N$, and secondly,
$(f,h) \to  \left|\sum_{i=1}^N (f(X_i)h(X_i)-\E f h) \right|$, the centred product process indexed by the classes $F$ and $H$.

We use chaining methods to obtain high probability upper bounds on the suprema of the two processes using a natural variation of Talagrand's $\gamma$-functionals.
\end{abstract}

\section{Introduction}
Empirical processes appear frequently in diverse branches of Mathematics, Statistics and Computer Science.

In its most standard form, an empirical process is indexed by a class of functions defined on a probability space $(\Omega,\mu)$. If $X_1,...,X_N$ are independent and distributed according to $\mu$, the centred empirical process indexed by $F$ is
$$
f \to \frac{1}{N}\sum_{i=1}^N f(X_i)-\E f, \ \ \ f \in F.
$$
One would like to obtain upper and lower bounds on
\begin{equation} \label{eq:empirical}
\sup_{f \in F} \left|\frac{1}{N}\sum_{i=1}^N f(X_i) -\E f \right|,
\end{equation}
either in high probability or in expectation. The hope is that the supremum \eqref{eq:empirical} may be controlled using some geometric features of the class $F$, similar to the ones used in the theory of gaussian processes (for more information on gaussian processes and their connection to the geometry of the underlying class see the books \cite{Dud-book} and \cite{Tal-book}).

\begin{Definition} \label{def:gauss-process}
Let $F \subset L_2$ and set $\{G_f : f \in F\}$ to be the centred canonical gaussian process indexed by $F$; that is, the gaussian process indexed by $F$ whose covariance structure is endowed by the inner product in $L_2$.
\end{Definition}

To avoid the (well understood) issue of measurability, set
$$
\E \sup_{f \in F} G_f = \sup \{ \E \sup_{h \in H} G_h \ : \ H \subset F, \ \ H \ {\rm is \ finite \ } \},
$$
and at times we shall write $\E\|G\|_F$ instead of $\E \sup_{f \in F} G_f$.
\vskip0.5cm
From certain perspectives, the standard empirical process is a more complex object that its gaussian counterpart. Roughly and somewhat inaccurately put, while a gaussian process $\{G_f : f \in F\}$ is governed by a single metric, $d(f,h)=(\E|G_f-G_h|^2)^{1/2}$, adequate control on an empirical process requires the use of a family of metrics. We will clarify this statement and develop it further, focusing on two natural generalizations of \eqref{eq:empirical}:

\begin{description}
\item{(1)} Let $\xi \in L_q$ for some $q > 2$ ($\xi$ need not be independent of $X$) and set $\xi_1,...,\xi_N$ to be independent copies of $\xi$. Consider
\begin{equation} \label{eq:multi-in-intro}
\sup_{f \in F} \left|\frac{1}{N}\sum_{i=1}^N \xi_i f(X_i) -\E \xi f \right|,
\end{equation}
which is the supremum of the multiplier process indexed by $F$ and associated with the multiplier $\xi$.

Note that unlike the standard notion of a multiplier process, here $(\xi_i)_{i=1}^N$ need not be independent of $(X_i)_{i=1}^N$.

Besides being a natural object from the theoretical point of view, the significance of multiplier processes may be seen in numerous applications. For example, multiplier processes play a central role in Statistics, when studying prediction and estimation problems (see, e.g. \cite{VW} and references therein, and also \cite{Men-LWC1,Men-LWC2}), but that is only the tip of the iceberg as far as applications go.

\item{(2)}  Let $F$ and $H$ be classes of functions defined on the probability space $(\Omega,\mu)$ and consider the supremum of the product process indexed by $F$ and $H$,
\begin{equation} \label{eq:prod-in-intro}
\sup_{f \in F, \ h \in H} \left|\frac{1}{N}\sum_{i=1}^N f(X_i)h(X_i) - \E f h \right|.
\end{equation}
Clearly, \eqref{eq:prod-in-intro} is a natural object when trying to analyze, for example, empirical correlation, or, what is arguably the most important process as far as applications are concerned, the quadratic empirical process
$$
f \to \frac{1}{N}\sum_{i=1}^N f^2(X_i) - \E f^2, \ \ \ f \in F.
$$
\end{description}

Although the product process may be viewed as the standard empirical process indexed by the product class $F \cdot H =\{f \cdot h : f \in F, h \in H\}$, and the multiplier process as the standard empirical process indexed by the class $\xi \cdot F=\{\xi f : f \in F\}$ (at least when $(\xi_i,X_i)_{i=1}^N$ are independent), the type of result one is looking for here is rather different. When studying the two processes, one would like to bound the supremum of \eqref{eq:multi-in-intro} and of \eqref{eq:prod-in-intro} using some geometric structures of the indexing classes $F$ and $F,H$ respectively, rather than the structures of $\xi \cdot F$ and $F \cdot H$, which, in most cases, are hard to handle.

\vskip0.5cm

Before we continue exploring the two problems, let us explain what is meant by ``some geometric structures of the indexing classes".

Motivated by {\it chaining} methods that have had tremendous impact on the theory of gaussian processes, the geometric parameters we shall focus on are `relatives' of Talagrand's $\gamma$ functionals.

\begin{Definition} \label{def:chaining}
Let $(T,d)$ be a metric space. An admissible sequence of $T$ is a collection of subsets, $T_s \subset T$, whose cardinality satisfies $|T_s| \leq 2^{2^s}$ for $s \geq 1$, and $|T_0|=1$. For $\alpha \geq 1$ and $s_0 \geq 0$ set
$$
\gamma_{s_0,\alpha}(T,d)=\inf \sup_{t \in T} \sum_{s \geq s_0} 2^{s/\alpha}d(t,T_s),
$$
where the infimum is taken with respect to all admissible sequences of $T$. When $s_0=0$ we shall write $\gamma_\alpha(T,d)$ instead of $\gamma_{0,\alpha}(T,d)$.
\end{Definition}
For more information on chaining methods we refer the reader to M. Talagrand's book \cite{Tal-book}, which contains an extensive and illuminating survey on the topic.

\vskip0.5cm
One of the most significant results in the theory of gaussian processes is based on chaining: that $\E \sup_{f \in F} G_f$ is characterized by $\gamma_2(F,L_2)$.
\begin{Theorem} \label{thm:MM}
There exist absolute constants $c$ and $C$ for which the following holds. Let $F \subset L_2$, $|F|>1$, and consider $\{G_f : f \in F\}$, the centred, canonical gaussian process indexed by $F$.
$$
c\gamma_2(F,L_2) \leq \E\|G\|_F \leq C \gamma_2(F,L_2).
$$
\end{Theorem}
The upper bound in Theorem \ref{thm:MM} is due to Fernique \cite{Fer} while the lower one is Talagrand's Majorizing Measures Theorem \cite{Tal-MM}. The proof of both parts may be found in \cite{Tal-book}. It should be noted that the original proofs of both results are based on the majorizing measures mechanism which preceded the modern generic chaining scheme.

\subsection{Chaining and the complexity parameter}
Chaining arises as a way of relating the supremum of a random process $\{Z_f : f \in F\}$ to the structure of the indexing set $F$. An upper estimate is obtained by combining individual tail bounds that ensure that if $f$ and $h$ are close in some sense, the probability that $Z_f$ is very different from $Z_h$ is small. For example, if $\{G_f : f \in F\}$ is the centred, canonical gaussian process indexed by $F \subset L_2$, and if $f,h \in F$, then
\begin{equation} \label{eq:gauss-tail}
Pr(|G_f-G_h| \geq u\|f-h\|_{L_2}) \leq 2\exp(-u^2/2),
\end{equation}
which is precisely the sort of tail estimate one would like to have (but unfortunately, analogous versions of \eqref{eq:gauss-tail} are not as simple for many interesting processes).

The {\it increment condition} \eqref{eq:gauss-tail} hints towards a fundamental fact: the supremum of a gaussian process $\{G_f : f \in F\}$ is determined by a single metric, endowed by the $L_2$ norm\footnote{In fact, Theorem \ref{thm:MM} implies a two-sided control using the $L_2$ metric, but since our focus is on obtaining upper bounds, we will focus only on that direction.}.
However, when it comes to empirical processes, the situation is rather different. To explain this substantial difference between empirical and gaussian processes it is convenient to use the notion of an Orlicz norm.
\begin{Definition} \label{def:psi-alpha}
Let $f$ be defined on the probability space $(\Omega,\mu)$. For $\alpha \geq 1$ set
$$
\|f\|_{\psi_\alpha}=\inf\{c>0 : \E \exp(|f/c|^\alpha) \leq 2\}
$$
and let $L_{\psi_\alpha}$ be the space of functions for which $\|f\|_{\psi_\alpha} < \infty$.
\end{Definition}
It is well known that $\|f\|_{\psi_\alpha}$ is equivalent to the smallest constant $c$ for which $Pr(|f|>t) \leq 2\exp(-t^\alpha/c^\alpha)$ for every $t>0$, and also to $\sup_{p \geq 1} \|f\|_{L_p}/p^{1/\alpha}$.
\vskip0.5cm

\begin{Definition} \label{def:subgauss}
A class $F$ is $L$-subgaussian if for every $f,h \in F \cup \{0\}$,
$$
\|f-h\|_{\psi_2} \leq L \|f-h\|_{L_2}.
$$
\end{Definition}

Using the moment characterization of the $\psi_2$ norm, it is evident that if $F$ is $L$-subgaussian then
$$
\|f-h\|_{L_p} \leq cL\sqrt{p}\|f-h\|_{L_2}
$$
for every $f,h \in F \cup \{0\}$ and for a suitable absolute constant $c$.

\vskip0.5cm

Observe that if $F \subset L_2$, the centred, canonical gaussian process $\{G_f : f \in F\}$ satisfies that $\|G_f-G_h\|_{L_2} = \|f-h\|_{L_2}$ and
$$
c_1\|G_f - G_h\|_{L_2} \leq \|G_f - G_h \|_{\psi_2} \leq  c_2L\|G_f -G_h\|_{L_2} =c_2L\|f-h\|_{L_2}
$$
for absolute constants $c_1$ and $c_2$. Unfortunately, such an equivalence need not be true for an empirical process. Indeed, if $Z_f = \sum_{i=1}^N (f(X_i)-\E f)$ then it is standard to verify that
$$
\|Z_f-Z_h\|_{\psi_2} \leq c_3\sqrt{N}\|f-h\|_{\psi_2},
$$
but there is no reason why $\|f-h\|_{\psi_2}$ should be equivalent to $\|f-h\|_{L_2}$ unless the class is subgaussian. And though one may show that there is an absolute constant $c_4$ for which
\begin{equation} \label{eq:subgaussian-emp-intro}
\E \sup_{f \in F} \left|\frac{1}{\sqrt{N}}\sum_{i=1}^N (f(X_i)-\E f) \right| \leq c_4 \gamma_2(F,\psi_2),
\end{equation}
the underlying metric in \eqref{eq:subgaussian-emp-intro} is the $\psi_2$ metric, which is simply too large to be of any use in most applications.

\vskip0.5cm

It turns out that this rather unsatisfactory upper bound may be improved by examining the chaining process more closely.

\vskip0.5cm

Let $\{Z_f : f \in F\}$ be a random process, set $s_0 \geq 0$ and consider an admissible sequence $(F_s)_{s \geq 0}$ and a collection of functions $\pi_s:F \to F_s$. Under mild assumptions (for example, that for every $f \in F$, $\pi_s f \to f$ in an appropriate sense), each $Z_f$ may be written as a `chain', which is simply a telescopic sum of the `links' $Z_{\pi_{s+1}f} -Z_{\pi_s f}$:
\begin{equation} \label{eq:chain}
Z_f = Z_{\pi_{s_0}f} + \sum_{s \geq s_0} \left(Z_{\pi_{s+1}f} - Z_{\pi_{s}f}\right).
\end{equation}
Obtaining uniform control over all chains requires balancing the tail estimates for each link that appears at the $s$-stage with the total number of links that are involved at that stage. And, since there are at most
$$
|F_s| \cdot |F_{s+1}| \leq 2^{2^s} \cdot 2^{2^{s+1}} \leq 2^{2^{s+2}}
$$
links at the $s$-stage, it suffices to find levels $t(f,s)$ for which
$$
Pr(|Z_{\pi_{s+1}f} - Z_{\pi_s f} | \geq t(f,s)) \leq \exp(-2^{s+2})
$$
and obtain a similar estimate for the `starting points' of each chain, $Z_{\pi_{s_0}f}$, to ensure the required uniform control over all possible chains. Such a uniform control results in a high probability upper estimate on $\sup_{f \in F} Z_f$.

A trivial yet crucial observation is that by Chebychev's inequality, a possible choice of $t(f,s)$ is
$$
t(f,s) \sim \|Z_{\pi_{s+1}f} - Z_{\pi_s f}\|_{L_{p}} \leq \|Z_{\pi_{s+1}f} -Z_f \|_{L_p} + \|Z_{\pi_s f} -Z_{f}\|_{L_p}
$$
for $p \sim 2^s$, where here and throughout the article we write $a \sim b$ if there are absolute constants $c_1$ and $c_2$ for which $c_1 a \leq b \leq c_2 a$.

Thus, an obvious alternative to the $\gamma$-functionals is
\begin{equation} \label{eq:comp-latala}
(*)=\inf \left(\sup_{f \in F} \sum_{s \geq s_0} \|Z_f-Z_{\pi_s f}\|_{L_{u^22^s}} + \|Z_{\pi_{s_0}f}\|_{L_{u^22^{s_0}}}\right),
\end{equation}
where the infimum is taken with respect to all admissible sequences of $F$ and for $\pi_s:F \to F_s$ which is the nearest point map with respect to the norm $\| \ \|_{L_{u^22^s}}$.

It immediately follows from the decomposition to chains in \eqref{eq:chain} that for $u \geq 4$, with probability at least $1-2\exp(-c_0u^22^{s_0})$,
\begin{equation} \label{eq:intro-diamond}
\sup_{f \in F} |Z_f| \leq c_1(*),
\end{equation}
where $c_0$ and $c_1$ are absolute constants.

\vskip0.5cm

The idea of using a complexity parameter that takes into account all the $L_p$ structures endowed by the process has been introduced in \cite{MenPao} and independently by R. Lata{\l}a  (see, for example, \cite{LatToc} and  Exercise 2.2.15 in \cite{Tal-book}).

\vskip0.5cm
Let us study two examples in which the way the complexity parameter in \eqref{eq:comp-latala} depends on $F$ takes a rather simple form.

\begin{description}
\item{$\bullet$} When considering the centred, canonical gaussian process indexed by $F$, the parameter in \eqref{eq:comp-latala} is not new. Indeed, recall that
$$
\|G_f-G_{\pi_s f}\|_{L_{u^22^s}} \sim u2^{2/s}\|f-\pi_s f\|_{L_2};
$$
hence, for an almost optimal admissible sequence \eqref{eq:comp-latala} becomes
\begin{equation} \label{eq:comp-outcome-intro}
\sim u \cdot \left(\sup_{f \in F} \sum_{s \geq s_0} 2^{s/2}\|f-{\pi_s f}\|_{L_2} + 2^{s_0/2}\|\pi_{s_0}f\|_{L_2}\right),
\end{equation}
and the relations between \eqref{eq:intro-diamond} and the upper estimate in Theorem \ref{thm:MM} via the $\gamma_2$ functional are clear.

\item{$\bullet$} Next, one may consider \eqref{eq:comp-latala} for the standard empirical process,
\begin{equation} \label{eq:empirical-process}
Z_f = \sum_{i=1}^N (f(X_i)-\E f).
\end{equation}
Applying Lata{\l}a's sharp bound on the moments of sums of independent random variables \cite{Lat-moments} (see also Theorem \ref{thm:latala}, below), one may show that for every $f,h \in F$ and every $p \geq 2$,
$$
\|Z_f - Z_h \|_{L_p} \leq c\sqrt{N} \sqrt{p} \cdot \sup_{1 \leq q \leq p} \frac{\|f-h\|_{L_q}}{\sqrt{q}}
$$
for an absolute constant $c$.
\end{description}

The last example leads to the introduction of the following norms and to a `graded version' of the $\gamma$ functionals.
\begin{Definition} \label{def:comp-norm}
For a random variable $Z$ and $p \geq 1$, set
$$
\|Z\|_{(p)} = \sup_{1 \leq q \leq p} \frac{\|Z\|_{L_q}}{\sqrt{q}}.
$$
\end{Definition}
Thus, if $Z_f$ is the empirical process from \eqref{eq:empirical-process}, it follows that
$$
\|Z_f-Z_h\|_{L_p} \leq c\sqrt{N} \sqrt{p} \|f-h\|_{(p)}.
$$
\begin{Definition} \label{def:comp-parameter}
Given a class of functions $F$, $u \geq 1$ and $s_0 \geq 0$, put
\begin{equation} \label{eq:comp-gamma-s-0}
{\Lambda}_{s_0,u}(F) = \inf \sup_{f \in F} \sum_{s \geq s_0} 2^{s/2} \|f - \pi_s f\|_{(u^22^s)},
\end{equation}
where the infimum is taken with respect to all admissible sequences $(F_s)_{s \geq 0}$, and $\pi_s f$ is the nearest point in $F_s$ to $f$ with respect to the $(u^22^s)$ norm.

Also, let
$$
\tilde{\Lambda}_{s_0,u}(F)={\Lambda}_{s_0,u}(F)+2^{s_0/2}\sup_{f \in F} \|\pi_{s_0}f\|_{(u^22^{s_0})}.
$$
\end{Definition}
Again, a straightforward chaining argument shows that for every $s_0 \geq 0$ and $u \geq 4$, with probability at least $1-2\exp(-c_0u^22^{s_0})$,
\begin{equation} \label{eq:simple-chaining-intro}
\sup_{f \in F} \left|\frac{1}{\sqrt{N}}\sum_{i=1}^N \left(f(X_i) - \E f\right) \right| \leq c_1 u \tilde{\Lambda}_{s_0,u}(F).
\end{equation}

Observe that if $F$ happens to be $L$-subgaussian and $u \geq 2$ then $\Lambda_{s_0,u}(F)$ is equivalent to $\gamma_{s_0,2}(F,L_2)$. Indeed, by the moment characterization of the $\psi_2$ norm, there is an absolute constant $c$ for which, for every $p \geq 2$, $\|f\|_{L_p} \leq c\sqrt{p}\|f\|_{\psi_2}$. Therefore,
$$
\frac{1}{\sqrt{2}}\|f-h\|_{L_2} \leq \|f-h\|_{(p)} \leq c\|f-h\|_{\psi_2} \leq cL\|f-h\|_{L_2},
$$
implying that
$$
\frac{1}{\sqrt{2}}\gamma_{s_0,2}(F,L_2) \leq \Lambda_{s_0,u}(F) \leq c\gamma_{s_0,2}(F,\psi_2) \leq cL \gamma_{s_0,2}(F,L_2).
$$
Moreover, if $\{G_f : f \in F\}$ is the centred, canonical gaussian process indexed by $F$, and $|F|>1$ then
by the Majorizing Measures Theorem $\gamma_2(F,L_2) \leq c_1 \E\|G\|_F$; since $\gamma_{s_0,2}(F,L_2) \leq \gamma_2(F,L_2)$, it follows that for $u \geq 1$,
$$
\tilde{\Lambda}_{s_0,u} (F)\leq c_2 L\left(\E\|G\|_F +2^{s_0/2}\sup_{f \in F} \|f\|_{L_2}\right).
$$

This leads to the next corollary, which may also be established directly, using a standard chaining argument.

\begin{Corollary} \label{cor:subgaussian-emp-intro}
Let $F$ be an $L$-subgaussian class. Then for every $u \geq 4$ and $s_0 \geq 0$, with probability at least $1-2\exp(-c_02^{s_0}u^2)$,
$$
\sup_{f \in F} \left|\frac{1}{\sqrt{N}}\sum_{i=1}^N \left(f(X_i)-\E f\right)\right| \leq c_1Lu \left(\E\|G\|_F + 2^{s_0/2}\sup_{f \in F}\|f\|_{L_2}\right).
$$

\end{Corollary}
As an example, let $d_2(F) = \sup_{f \in F} \|f\|_{L_2}$
and set $s_0 \geq 0$ be the largest integer for which
$$
\E\|G\|_F \geq 2^{s_0/2}d_2(F);
$$
if no such integer exists, set $s_0=0$ and note that $\E\|G\|_F \geq c d_2(F)$. Hence, with probability at least
$$
1-2\exp\left(-c_2u^2 \left(\frac{\E\|G\|_F}{d_2(F)}\right)^2\right),
$$
$$
\sup_{f \in F} \left|\frac{1}{\sqrt{N}}\sum_{i=1}^N \left(f(X_i)-\E f\right)\right| \leq c_3Lu \E\|G\|_F.
$$
\vskip1cm
However, and unlike subgaussian examples, there are natural examples in which ${\Lambda}_{s_0,u}(F)$ may be significantly smaller than its $\psi_2$ counterpart. This should not come as a surprise, as both $\| \ \|_{L_{2^s}}/2^{s/2}$ and $\| \ \|_{(2^s)}$ are `local' versions of the $\psi_2$ metric: they measure the subgaussian behaviour of the functions involved, but only up to a fixed level, rather than at every level.

\vskip0.5cm
\noindent{\bf Example.} Let $T \subset \R^n$ and consider $F_T = \left\{\inr{t,\cdot} : t \in T\right\}$, the class of linear functionals indexed by $T$. Let $Y=(y_1,...,y_n)$ be a random vector with independent, standard exponential random variables as coordinates and set $\mu$ to be the underlying measure endowed on $\R^n$ by $Y$.

Clearly, if $T$ contains any one of the coordinate directions $\{e_1,...,e_n\}$, then $F_T$ is not a subset of $L_{\psi_2}$ and $\gamma_2(F,\psi_2)$ is not even well defined.

On the other hand, in \cite{GluKwa}, Gluskin and Kwapien showed that for every $t \in \R^n$ and $p \geq 1$,
$$
\|\inr{t,Y}\|_{L_p} \sim p\|t\|_{\ell_\infty^n} + \sqrt{p} \|t\|_{\ell_2^n}.
$$
Therefore,
$$
\|\inr{t,Y}\|_{(p)} = \sup_{1 \leq q \leq p} \frac{\|\inr{t,Y}\|_{L_q}}{\sqrt{q}} \sim \sqrt{p} \|t\|_{\ell_\infty^n} + \|t\|_{\ell_2^n};
$$
$\|\inr{t,Y}\|_{(u^22^s)} \sim u 2^{s/2} \|t\|_{\ell_\infty^n} + \|t\|_{\ell_2^n}$; and
\begin{align} \label{eq:lambda-log-concave}
{\Lambda}_{s_0,u}(F) \sim & \inf \sup_{t \in T} \left(\sum_{s \geq s_0} u2^s\|t-\pi_s t\|_{\ell_\infty^n} + 2^{s/2} \|t-\pi_s t\|_{\ell_2^n}\right) \nonumber
\\
\sim & u\gamma_{s_0,1}(T,\ell_\infty^n) + \gamma_{s_0,2}(T,\ell_2^n).
\end{align}
Talagrand showed in \cite{Tal-AM} (see also \cite{Tal-book}) that \eqref{eq:lambda-log-concave} has a geometric interpretation:
$$
\gamma_{1}(T,\ell_\infty^n) + \gamma_{2}(T,\ell_2^n) \sim \E \sup_{t \in T} \inr{t,Y},
$$
which is the mean-width of $T$ relative to the random vector $Y$. Thus, and in contrast to the $\psi_2$-based parameter, there is an absolute constant $c$ for which, for every $T \subset \R^n$ and every $u \geq 1$,
\begin{equation} \label{eq:tal-exp-control}
\tilde{\Lambda}_{0,u}(F_T) \leq c u \E \sup_{t \in T} \inr{t,Y}.
\end{equation}

\subsection{The main results}
Up to now, we have only examined the standard empirical process. Unfortunately, the simple chaining argument used in \eqref{eq:simple-chaining-intro} to control that process is rather useless when it comes to dealing with multiplier processes or with product processes. We will show in what follows how the suprema of these two types of processes may be bounded from above in terms of the $\Lambda$-functionals.

\vskip0.5cm
Let us begin by formulating the estimate for multiplier processes.

\begin{Theorem} \label{thm:multiplier-intro}
For $q>2$, there are constants $c_0$, $c_1,c_2$ and $c_3$ that depend only on $q$, for which the following holds. Let $\xi \in L_q$ and set $\xi_1,...,\xi_N$ to be independent copies of $\xi$. Fix an integer $s_0 \geq 0$ and $w,u>c_0$. Then, with probability at least
$$
1-c_1w^{-q} N^{-((q/2)-1)}\log^{q} N-2\exp(-c_2u^2 2^{s_0}),
$$
$$
\sup_{f \in F} \left|\frac{1}{\sqrt{N}} \sum_{i=1}^N \left(\xi_i f(X_i) - \E \xi f\right) \right| \leq c_3wu\|\xi\|_{L_q}\tilde{\Lambda}_{s_0,u}(F).
$$
\end{Theorem}

One simple outcome of Theorem \ref{thm:multiplier-intro} is when the class $F$ happens to be $L$-subgaussian.
\begin{Corollary} \label{cor:multi-subgaussian-intro}
Recall that $d_2(F) = \sup_{f \in F} \|f\|_{L_2}$ and set $s_0 \geq 0$ to satisfy that
$$
\E\|G\|_F \sim 2^{s_0/2}d_2(F).
$$
 As noted previously, since $F$ is $L$-subgaussian,
$$
\tilde{\Lambda}_{s_0,u}(F) \leq c_1L \E\|G\|_F.
$$
Thus, it follows from Theorem \ref{thm:multiplier-intro} that with probability at least
$$
1-c_2(q)w^{-q} N^{-((q/2)-1)}\log^{q} N-2\exp(-c_3u^2 (\E\|G\|_F/d_2(F))^2),
$$
\begin{equation} \label{eq:multi-subgaussian-intro}
\sup_{f \in F} \left|\frac{1}{\sqrt{N}}\sum_{i=1}^N \left(\xi_i f(X_i) - \E \xi f\right) \right| \leq c_4(q) Lwu\|\xi\|_{L_q}\E\|G\|_F.
\end{equation}
\end{Corollary}
\vskip0.5cm

Turning to the question of product processes, recall the following fact from \cite{MenPao} (see also Theorem 9.3.1 in \cite{Tal-book}).

\begin{Theorem} \label{thm:quad-talagrand}
For every $q>4$ there exists a constant $c(q)$ that depends only on $q$ for which the following holds. Let $F \subset L_q$ be a class of functions on $(\Omega,\mu)$. Assume that for every $t>0$ and every $f,h \in F \cup \{0\}$,
$$
Pr(|f-h| \geq t) \leq 2\exp\left(-\min\left\{\frac{t^2}{d_2^2(f,h)},\frac{t}{d_1(f,h)}\right\}\right),
$$
where $d_1$ and $d_2$ are metrics on $F \cup \{0\}$. If $d_q(F)=\sup_{f \in F} \|f\|_{L_q}$ and $\gamma=\gamma_2(F,d_2)+\gamma_1(F,d_1)$, then
$$
\E \sup_{f \in F} \left|\frac{1}{N} \sum_{i=1}^N f^2(X_i) - \E f^2\right| \leq c_q\left(d_q(F) \frac{\gamma}{\sqrt{N}}+\frac{\gamma^2}{N}\right).
$$
\end{Theorem}

Theorem \ref{thm:quad-talagrand} is demonstrated by integrating a high-probability bound. However, the probability estimate established in \cite{MenPao} is far from optimal. Recently, Dirksen \cite{Dirk} obtained the optimal probability estimate under the assumption that $F \subset L_{\psi_2}$, improving earlier results from \cite{MPT,Men-psi1}, and a few months later, Bednorz \cite{Bed} gave a different proof of the same fact. The following formulation is from \cite{Bed}.

\begin{Theorem} \label{thm:Bednorz}
There exist absolute constants $c_1$ and $c_2$ for which the following holds. Let $F$ be a class of functions on $(\Omega,\mu)$ and set $\gamma=\gamma_2(F,\psi_2)$ and $d_{\psi_2}(F)=\sup_{f \in F} \|f\|_{\psi_2}$. For every $u >0$, with probability at least $1-2\exp(-c_1\min\{\sqrt{N}u,u^2\})$,
\begin{equation} \label{eq:Bed}
\sup_{f \in F} \left|\frac{1}{N} \sum_{i=1}^N f^2(X_i) - \E f^2 \right| \leq c_2\left(d_{\psi_2}(F) \frac{\gamma}{\sqrt{N}}+\frac{\gamma^2}{N}+ u \frac{d^2_{\psi_2}(F)}{\sqrt{N}} \right).
\end{equation}
\end{Theorem}
The probability estimate in Theorem \ref{thm:Bednorz} is indeed optimal, as may be seen by setting $t=ud_{\psi_2}^2(F)/\sqrt{N}$ in Bernstein's inequality applied to $f^2$.

\vskip0.5cm

In comparison, below is our estimate on the supremum of a product process, and in particular, on the supremum of the quadratic process.

\begin{Theorem} \label{thm:quadratic-process-upper-intro}
There exists an absolute constant $c_0$ and for every $q>4$ there exists a constant $c_1(q)$ that depends only on $q$ for which the following holds. Let $F$ and $H$ be classes of functions on $(\Omega,\mu)$, set $u \geq \max\{8,\sqrt{q}\}$ and consider an integer $s_0 \geq 0$. Then, with probability at least $1-2\exp(-c_0u^22^{s_0})$, for every $f \in F$ and $h \in H$,
\begin{align*}
& \left|\sum_{i=1}^N f(X_i)h(X_i) - \E fh \right|
\\
\leq & c_1(q) \left( u^2 \tilde{\Lambda}_{s_0,u}(H)\tilde{\Lambda}_{s_0,u}(F) + u\sqrt{N}\left(d_q(F)\tilde{\Lambda}_{s_0,u}(H)+d_q(H)\tilde{\Lambda}_{s_0,u}(F)\right)\right).
\end{align*}
In particular, if $F=H$, then with probability at least $1-2\exp(-c_0u^22^{s_0})$,
\begin{equation} \label{eq:quad-intro-main}
\left|\sum_{i=1}^N (f^2(X_i)-\E f^2) \right| \leq c_1(q) \left(u^2 \tilde{\Lambda}_{s_0,u}^2(F)+u \sqrt{N} d_q(F)\tilde{\Lambda}_{s_0,u}(F)\right).
\end{equation}
\end{Theorem}

\vskip0.5cm

Let us present some of the outcomes of Theorem \ref{thm:quadratic-process-upper-intro} for the quadratic process.
\begin{description}
\item{$\bullet$} Since $\| \ \|_{(p)} \leq c \| \ \|_{\psi_2}$, it is evident that $\Lambda_{s_0,u}(F) \leq c \gamma_2(F,\psi_2)$. Thus,
    $$
    \tilde{\Lambda}_{s_0,u}(F) \leq \Lambda_{s_0,u}(F)+2^{s_0/2}d_{\psi_2}(F) \leq c_1\gamma_2(F,\psi_2)
    $$
if $|F|>1$ and one sets $2^{s_0} \sim (\gamma_2(F,\psi_2)/d_{\psi_2}(F))^2$. Also, since $\|f\|_{\psi_2} \sim \sup_{q \geq 2} \|f\|_{L_q}/\sqrt{q}$, it is evident that $d_q(F) \leq c_2\sqrt{q}d_{\psi_2}(F)$, and \eqref{eq:quad-intro-main} recovers Theorem \ref{thm:Bednorz}.
\item{$\bullet$} If the indexing class is $L$-subgaussian, then by applying the same argument as in Corollary \ref{cor:multi-subgaussian-intro} for the choice $2^{s_0} \sim (\E\|G\|_F/d_2(F))^2$, it is evident that with probability at least $1-2\exp(-c_0u^2(\E\|G\|_F/d_2(F))^2)$,
\begin{equation} \label{eq:quad-subgaussian-intro}
\left|\sum_{i=1}^N (f^2(X_i)-\E f^2) \right| \leq c_1 L^2\left(u^2 \left(\E\|G\|_F\right)^2+u \sqrt{N} d_2(F)\E\|G\|_F\right).
\end{equation}
\end{description}
\vskip0.5cm

The proofs of Theorem \ref{thm:multiplier-intro} and Theorem \ref{thm:quadratic-process-upper-intro} are based on symmetrization, which has been one of the most influential tools in empirical processes theory. The most well-known symmetrization inequalities for empirical processes are the celebrated Gin\'{e}-Zinn inequalities \cite{GZ}, but we will use an earlier, ``in-probability" version of those inequalities (see, e.g., \cite{VW}).

\begin{Theorem} \label{thm:GZ}
Let $(Z_f(i))_{i=1}^N$ be independent copies of a mean-zero stochastic process $\{Z_f : f \in F\}$, and for every $1 \leq i \leq N$, set $y_f(i):F \to \R$ to be arbitrary functions. Let $(\eps_i)_{i=1}^N$ be independent, symmetric, $\{-1,1\}$-valued random variables that are independent of $(Z_f(i))_{i=1}^N$. Then, for every $x>0$
\begin{align*}
\left(1-\frac{4N}{x^2} \sup_{f \in F}{\rm var}\left(Z_f\right)
\right) & Pr \left( \sup_{f \in F} \left|\sum_{i=1}^N Z_f(i)\right| > x
\right)
\\ \leq & 2Pr \left(\sup_{f \in F} \left|\sum_{i=1}^N
\eps_i\left(Z_f(i)-y_f(i)\right)\right| > \frac{x}{4} \right).
\end{align*}
\end{Theorem}

A symmetrization result may be derived for the standard empirical processes by setting $Z_f = f(X) - \E f$ and $y_f(i)=\E f$;  for the multiplier process by setting $Z_f = \xi f(X) - \E \xi f$ and $y_f(i)=-\E \xi f$; and for the product process by setting $Z_{f,h} = f(X)h(X)-\E fh$ and $y_{f,h}(i)=-\E fh$.

\vskip0.5cm

Thanks to Theorem \ref{thm:GZ}, one may prove Theorem \ref{thm:multiplier-intro} via a high-probability upper bound on the supremum of the Bernoulli process
$$
\sup_{v \in V} \left|\sum_{i=1}^N \eps_i z_i v_i \right|,
$$
where the set $V$ is a typical coordinate projection of the class $F$; that is, for $\sigma=(X_1,...,X_N)$,
$$
V=P_\sigma F= \left\{(f(X_i))_{i=1}^N : f \in F \right\},
$$
and $z \in \R^N$ is a typical realization of the random vector $(\xi_i)_{i=1}^N$.

In a similar fashion, Theorem \ref{thm:quadratic-process-upper-intro} follows from a high-probability upper bound on
$$
\sup_{v \in V, w \in W} \left|\sum_{i=1}^N \eps_i w_i v_i \right|,
$$
for typical coordinate projections $V=P_\sigma F$ and $W=P_\sigma H$.

\vskip0.5cm
This observation dictates the structure of the article. We will first study
$$
\sup_{v \in V} \left|\sum_{i=1}^N \eps_i z_i v_i \right| \ \ {\rm and} \ \ \sup_{v \in V, w \in W} \left|\sum_{i=1}^N \eps_i w_i v_i \right|
$$
for fixed sets $V,W \subset \R^N$ and $z \in \R^N$, focusing on the way in which the geometry of the indexing sets is manifested in the main complexity parameter that appears in the upper bounds. We will then explore this (deterministic) complexity parameter and show that for a typical coordinate projection of each indexing set, it is, in fact, $\tilde{\Lambda}_{s_0,u}(F)$.

\vskip0.5cm

Finally, a word about notation. Throughout, $c_0,c_1,...$ denote absolute constants. Their value may change from line to line. $c(q)$ or $c_q$ are constants that depend only on the parameter $q$, and $a \lesssim_q b$ means that $a \leq c(q) b$.

For $1 \leq p \leq \infty$, denote by $\ell_p^m$ the space $\R^m$ endowed with the $\ell_p$ norm. And, for $(x_i)_{i=1}^m \in \R^m$, let $(x_i^*)_{i=1}^m$ be the non-increasing rearrangement of $(|x_i|)_{i=1}^m$.

\section{Chaining and Bernoulli processes} \label{sec:chain-and-Bernoulli}
As we noted earlier, our method of analysis consists of two main components. First, gathering accurate information on the structure of a typical coordinate projection of $F$ (and in the case of the product process, of a typical coordinate projection of $H$ as well); and second, for a typical $\sigma=(X_1,...,X_N)$ and $(\xi_i)_{i=1}^N$, analyzing the suprema of the conditioned Bernoulli processes
\begin{equation*}
f \to \left|\sum_{i=1}^N \eps_i \xi_i f(X_i)\right| \ \Big{|} \sigma, (\xi_i)_{i=1}^N
\end{equation*}
and
\begin{equation*}
(f,h) \to \left|\sum_{i=1}^N \eps_i f(X_i)h(X_i)\right| \ \Big{|} \sigma.
\end{equation*}

This section focuses on the latter: if $V, W \subset \R^N$ and $z =(z_i)_{i=1}^N \in \R^N$, we are interested in high probability upper bounds on
$$
\sup_{v \in V} \left|\sum_{i=1}^N \eps_i z_i v_i \right|
$$
and
$$
\sup_{v \in V, w \in W} \left|\sum_{i=1}^N \eps_i v_i w_i \right|.
$$

The upper bounds in both cases are based on a chaining argument combined with H\"{o}ffding's inequality. The key observation is that if $x \in \R^N$ and $(x_i^*)_{i=1}^N$ is the nonincreasing rearrangement of $(|x_i|)_{i=1}^N$, then for every $t>0$ and every fixed $1 \leq k < N$,
$$
Pr\left( \left|\sum_{i=1}^N \eps_i x_i \right| \geq \sum_{i=1}^k x_i^* + t \left(\sum_{i=k+1}^N (x_i^*)^2\right)^{1/2} \right) \leq 2\exp(-t^2/2).
$$
Moreover, by \cite{MontSmith}, this estimate is optimal when $k \sim t^2$.

Therefore, the effect $x$ has on the Bernoulli process depends on
$$
\sum_{i=1}^k x_i^* \ \ {\rm and} \ \ \left(\sum_{i=k+1}^N (x_i^*)^2\right)^{1/2},
$$
that is, the $\ell_1^N$ norm of the largest $k$ coordinates of $(|x_i|)_{i=1}^N$ and the $\ell_2^N$ norm of the smallest $N-k$ coordinates.

When dealing with products, as we have to, $x_i=w_i v_i$ and the decomposition requires additional care. Let $I$ be the union of the sets of the $k$ largest coordinates of $(|w_i|)_{i=1}^N$ and the $k$ largest coordinates of $(|v_i|)_{i=1}^N$; thus $|I| \leq 2k$. Fix $r>1$, set $r^\prime$ to be its conjugate index (i.e., $1/r+1/r^\prime=1$) and observe that with probability at least $1-2\exp(-t^2/2)$,
\begin{align} \label{eq:prod-monotone-basic}
& \left|\sum_{i=1}^N \eps_i x_i \right| \leq  \sum_{i \in I} x_i^* + t \left(\sum_{i \in I^c} (x_i^*)^2\right)^{1/2}
\\
\leq & 2\left(\sum_{i=1}^k (w_i^*)^2\right)^{1/2} \left(\sum_{i=1}^k (v_i^*)^2\right)^{1/2}
+ t \left(\sum_{i=k+1}^N (w_i^*)^{2r}\right)^{1/2r} \left(\sum_{i=k+1}^N (v_i^*)^{2r^\prime}\right)^{1/2r^\prime}. \nonumber
\end{align}

As \eqref{eq:prod-monotone-basic} is meant to play a part in a chaining argument, the bound should hold uniformly for all the links that appear at the $s$-stage; hence, a likely choice in \eqref{eq:prod-monotone-basic} at the $s$-stage is $t \sim 2^{s/2}$.

\subsection{The Bernoulli multiplier process} \label{sec:Bernoulli-multi}
Let us present a chaining process aimed at bounding
$$
\sup_{v \in V} \left|\sum_{i=1}^N \eps_i z_i v_i\right|,
$$
for $z \in \R^N$ and $V \subset \R^N$.

Let $(V_s)_{s \geq 0}$ be an admissible sequence of $V$, and note that under very mild assumptions (for example, that for every $v \in V$, $\pi_s v \to v$ in an appropriate sense as $s$ tends to infinity), for every $s_0 \geq 0$,
$$
v=\sum_{s \geq s_0} (\pi_{s+1}v - \pi_s v) + \pi_{s_0} v.
$$
Set $\Delta_s v = \pi_{s+1}v -\pi_s v$; thus
\begin{equation*}
\left|\sum_{i=1}^N \eps_i z_i v_i \right| \leq  \left|\sum_{s \geq s_0} \sum_{i=1}^N \eps_i z_i \cdot (\Delta_s v)_i \right|+\left|\sum_{i=1}^N \eps_i z_i \cdot (\pi_{s_0}v)_i \right|.
\end{equation*}
Let $r,r^\prime>1$ be conjugate indices, for every  $s \geq s_0$ set an integer $j_s \geq 1$ to be determined later and assume that $(j_s)_{s \geq s_0}$ is non-decreasing in $s$. Finally, let $I=I_{v,s}$ be the union of the $j_s-1$ largest coordinates of $(|\Delta_s v|_i)_{i=1}^N$ and the $j_s-1$ largest coordinates of $(|z_i|)_{i=1}^N$. Applying \eqref{eq:prod-monotone-basic} for $t \geq 4$, it is evident that with probability at least $1-2\exp(-t^2 2^s/2)$ 
\begin{align*}
& \left|\sum_{i=1}^N \eps_i z_i \cdot (\Delta_s v)_i \right| \leq
\\
\leq & 2 \|z\|_{\ell_2^N} \left(\sum_{i < j_s} ((\Delta_s v)^2)_i^*\right)^{1/2} + t 2^{s/2} \left(\sum_{i \geq j_s} (z_i^*)^{2r}\right)^{1/2r} \cdot \left(\sum_{i \geq j_s} ((\Delta_s v)_i^*)^{2r^\prime}\right)^{1/2r^\prime}.
\end{align*}
Repeating this argument for the vectors $\pi_{s_0}v$ and summing the probabilities, it follows that with probability at least $1-2\exp(-ct^22^{s_0})$, for every $v \in V$,
\begin{align} \label{eq:comp-bernoulli-multi}
& \left|\sum_{i=1}^N \eps_i z_i v_i\right| \leq
\\
& 2\|z\|_{\ell_2^N} \cdot \left(\sum_{s \geq s_0} \left(\sum_{i < j_s} ((\Delta_s v)_i^*)^2\right)^{1/2} + \left(\sum_{i < j_{s_0}} ((\pi_{s_0} v)_i^*)^2\right)^{1/2}\right) \nonumber
\\
+& t \left(\sum_{s \geq s_0} 2^{s/2} \left(\sum_{i \geq j_s} (z_i^*)^{2r}\right)^{1/2r} \cdot \left(\sum_{i \geq j_s} ((\Delta_s v)_i^*)^{2r^\prime}\right)^{1/2r^\prime}\right)
\nonumber
\\
+ & t 2^{s_0/2} \left(\sum_{i \geq j_{s_0}} (z_i^*)^{2r}\right)^{1/2r} \cdot \left(\sum_{i \geq j_{s_0}} ((\pi_{s_0} v)_i^*)^{2r^\prime}\right)^{1/2r^\prime}.
\nonumber
\end{align}
\vskip0.5cm

Motivated by this chaining argument, consider the following structural assumption on $V$:
\begin{Assumption} \label{ass:decomp-of-a-set}
Let $V \subset \R^N$, set $p \geq 1$, fix an integer $s_0$ and let $(V_s)_{s \geq 0}$ be an admissible sequence of $V$.
\begin{description}
\item{$\bullet$} Let $(\| \ \|_{[s]})$ be a family of semi-norms on $\R^N$, and set $\| \ \|$ to be an additional semi-norm on $\R^N$.
\item{$\bullet$} Let $(j_s)_{s \geq s_0}$ be a non-decreasing sequence of integers, and for every $s \geq s_0$, $1 \leq j_s \leq N+1$.
\end{description}
Assume that for every $s \geq s_0$ and every $v \in V$,
\begin{description}
\item{$\bullet$} $\left(\sum_{i < j_s} \left((\Delta_s v)_i^*\right)^2\right)^{1/2} \leq \|\Delta_s v\|_{[s]}$, \ \ $\left(\sum_{i \geq j_s} \left((\Delta_s v)_i^*\right)^{2p}\right)^{1/2p} \leq \|\Delta_s v\|N^{1/2p}$,
\item{$\bullet$} $\left(\sum_{i < j_s} \left((\pi_{s} v)_i^*\right)^2\right)^{1/2} \leq \|\pi_{s} v\|_{[s]}$, \ \  $\left(\sum_{i \geq j_{s}} \left((\pi_{s} v)_i^*\right)^{2p}\right)^{1/2p} \leq \|\pi_{s} v\|N^{1/2p}$.
\end{description}
\end{Assumption}

Observe that if $V$ satisfies Assumption \ref{ass:decomp-of-a-set} for $p=r^\prime$,
then by \eqref{eq:comp-bernoulli-multi}, with probability at least $1-2\exp(-ct^22^{s_0})$, for every $v \in V$,
\begin{align*}
 \left|\sum_{i=1}^N \eps_i z_i v_i \right|
\leq & 2 \|z\|_{\ell_2^N} \left(\sum_{s \geq s_0} \|\Delta_s v\|_{[s]} + \|\pi_{s_0}v\|_{[s_0]}\right)
\\
+ & t \left(\sum_{s \geq s_0} 2^{s/2}\|\Delta_s v\| + 2^{s_0/2} \|\pi_{s_0} v\|\right) \cdot N^{1/2r^\prime} \left(\sum_{i \geq j_{s_0}} (z_i^*)^{2r}\right)^{1/2r}.
\end{align*}

With this in mind, set
$$
\Lambda(V)=\sup_{v \in V} \left(\sum_{s \geq s_0} \|\Delta_s v\|_{[s]} + \|\pi_{s_0}v\|_{[s_0]}\right),
$$
$$
\Theta(V)=\sup_{v \in V} \left(\sum_{s \geq s_0} 2^{s/2}\|\Delta_s v\| + 2^{s_0/2}\|\pi_{s_0}v\|\right)
$$
and
$$
d(V)= \sup_{v \in V} \|v\|.
$$
\begin{Corollary} \label{cor:multi-in-section}
Using the same notation as above, with probability at least $1-2\exp(-ct^22^{s_0})$, for every $v \in V$
\begin{equation} \label{eq:Bernoulli-multiplier-clean}
\left|\sum_{i=1}^N \eps_i z_i v_i \right| \leq 2\|z\|_{\ell_2^N}\Lambda(V) + t\Theta(V)N^{1/2r^\prime} \left(\sum_{i \geq j_{s_0}} (z_i^*)^{2r}\right)^{1/2r}.
\end{equation}
\end{Corollary}
\vskip0.5cm

Seemingly, there is plenty of freedom in the choices of $j_s$, $r$ and the admissible sequence of $V$. However, our main interest is when $z$ and $V$ are typical realizations of $(\xi_i)_{i=1}^N$ and $V=P_\sigma F$ respectively, and the natural choice of an admissible sequence of $V$ should be endowed by an admissible sequence of the underlying class $F$. Thus, one must have adequate control on all the `monotone sums'
\begin{equation} \label{eq:monotone-cond-1}
\left(\sum_{i \in I_s } (\Delta_s f)^2(X_i)\right)^{1/2}, \ \ \left(\sum_{i \in I_s^c} (\Delta_s f)^{2r^\prime}(X_i)\right)^{1/2r^\prime},
\end{equation}
where $\Delta_s f = \pi_{s+1}f -\pi_s f$, $I_s$ is the set of the $j_s-1$ largest coordinates of $\left(|\Delta_s f|(X_i)\right)_{i=1}^N$. In a similar fashion, one must be able to control
\begin{equation} \label{eq:monotone-cond-2}
\left(\sum_{i \in I_{s_0}} (\pi_{s_0} f)^2(X_i)\right)^{1/2}, \ \ \left(\sum_{i \in I_{s_0}^c} (\pi_{s_0} f)^{2r^\prime}(X_i)\right)^{1/2r^\prime}.
\end{equation}

Since $\{\Delta_s f : f \in F\}$ contains at most $2^{2^{s+2}}$ points that must be controlled  uniformly, the individual probability estimate that is required in \eqref{eq:monotone-cond-1} and in \eqref{eq:monotone-cond-2} is $\exp(-u^22^s)$ for a large enough $u$ (a choice of $u \geq 4$ will do). In what follows, we will show that this almost forces the choice of $j_s$.

In addition, obtaining sufficient control on $\left(\sum_{i = j_s}^N (\xi_i^*)^{2r} \right)^{1/2r}$ for every $j_s$ restricts the choice of $r$; it will depend on the $L_q$ space to which $\xi$ belongs.

\vskip0.5cm
The key estimate on
$$
\left(\sum_{i < j} (Z_i^*)^2\right)^{1/2} \ \ {\rm and} \ \ \left(\sum_{i \geq j} (Z_i^*)^r\right)^{1/r}
$$
for $N$ independent copies of a random variable $Z$ will be derived in Section \ref{sec:pre}.

\subsection{Bernoulli product processes}
Chaining for a Bernoulli product process is more involved than the one outlined above. However, the two share a common feature: structural assumptions on the indexing sets.

\vskip0.5cm

Let $V,W \subset \R^N$ and assume that both $V$ and $W$ satisfy the following for $p=1$ and $p=2$.

\begin{Assumption} \label{ass:quad-strucutre}
Let $V \subset \R^N$, fix an integer $s_0$ and let $(V_s)_{s \geq 0}$ be an admissible sequence of $V$.
\begin{description}
\item{$\bullet$} Let $(\| \ \|_{[s]})$ be a family of semi-norms on $\R^N$ and set $\| \ \|$ to be an additional semi-norm on $\R^N$.
\item{$\bullet$} Let $(j_s)_{s \geq s_0}$ be an non-decreasing sequence of integers and for every $s \geq s_0$, $1 \leq j_s \leq N+1$.
\item{$\bullet$} Let $(\nu_s)_{s \geq s_0}$ be a sequence of nonnegative numbers.
\end{description}
Assume that for every $s \geq s_0$ and every $v \in V$,
\begin{description}
\item{$\bullet$} $\left(\sum_{i < j_s} \left((\Delta_s v)_i^*\right)^2\right)^{1/2} \leq \|\Delta_s v\|_{[s]}$, \ \ $\left(\sum_{i \geq j_{s-1}} \left((\Delta_s v)_i^*\right)^{2p}\right)^{1/2p} \leq \|\Delta_s v\|N^{1/2p}$,
\item{$\bullet$} $\left(\sum_{i < j_s} \left((\pi_{s} v)_i^*\right)^2\right)^{1/2} \leq \|\pi_{s} v\|_{[s]}$, \ \  $\left(\sum_{i  \geq j_{s-1}} \left((\pi_{s} v)_i^*\right)^{2p}\right)^{1/2p} \leq \|\pi_{s} v\|N^{1/2p}$.
\item{$\bullet$} $\left(\sum_{i  = j_{s-1}}^{j_s-1} \left((\pi_{s} v)_i^*\right)^{2}\right)^{1/2} \leq d(V) \nu_s$, where we set $j_{s_0-1}=j_{s_0}$ and recall that $d(V)=\sup_{v \in V} \|v\|$. 
\end{description}
\end{Assumption}

\begin{Remark}
There is a slight overlap between the sets of large and small coordinates: one set contains the $j_s-1$ largest coordinates, while the other contains the $N-j_{s-1}+1$ smallest coordinates -- rather than the more natural set, consisting of the $N-j_s+1$ smallest coordinates. The reason for this overlap is a minor technicality that will be used in the proof of Lemma \ref{lemma:ell-2-N-radius}.
\end{Remark}

\begin{Lemma} \label{lemma:ell-2-N-radius}
Consider $V \subset \R^N$ that satisfies Assumption \ref{ass:quad-strucutre} for $p=1$. Set $s_1 > s_0$ to be the first integer for which $j_s = N+1$, and if no such integer exists, set $s_1=s_0+1$. Then
$$
\sup_{v \in V} \left(\sum_{i=1}^N v_i^2 \right)^{1/2} \leq \Lambda(V)+ d(V) \cdot \left(\sqrt{N}+\sum_{s=s_0}^{s_1-1} \nu_s\right) .
$$
In particular, if $\sum_{s=s_0}^{s_1-1} \nu_s \leq c\sqrt{N}$ then
$$
\sup_{v \in V} \left(\sum_{i=1}^N v_i^2 \right)^{1/2} \leq \Lambda(V)+ (c+1)d(V)\sqrt{N}.
$$
\end{Lemma}

\proof  We will only consider the case $s_1>s_0+1$, as the proof of the case $s_1=s_0+1$ is simpler, and is actually contained in the proof of the former.
\vskip0.5cm
Fix $\pi_{s_1}v $ and let $I_{s_1-1}$ be the set of the $j_{s_1-1}-1$ largest coordinates of  $(|(\pi_{s_1}v)_i|)_{i=1}^N$. Therefore,
\begin{align*}
& \left(\sum_{i=1}^N (\pi_{s_1}v)_i^2\right)^{1/2}
\\
\leq & \left(\sum_{i \in I_{s_1-1}^c} (\pi_{s_1}v)_i^2\right)^{1/2}+\left(\sum_{i \in I_{s_1-1}} (\Delta_{s_1-1}v)_i^2\right)^{1/2}
+ \max_{|J|=j_{s_1-1}} \left(\sum_{j \in J} (\pi_{s_1-1}v)_j^2\right)^{1/2}.
\end{align*}
Applying Assumption \ref{ass:quad-strucutre} for $p=1$,
$$
\left(\sum_{i \in I_{s_1-1}^c} (\pi_{s_1}v)_i^2\right)^{1/2} = \left(\sum_{i \geq j_{{s_1}-1}} ((\pi_{s_1}v)_i^*)^2\right)^{1/2} \leq \|\pi_{s_1}v\| N^{1/2},
$$
and 
$$
\left(\sum_{i \in I_{s_1-1}} (\Delta_{s_1-1}v)_i^2\right)^{1/2} \leq  \|\Delta_{s_1-1}v\|_{[s_1-1]}.
$$
Repeating this argument for $s_0 < s <s_1$, it follows that for every $v \in V$, and for the choice of $I_s$ as the set of the $j_s-1$ largest coordinates of $(|\pi_{s+1}v|_i)_{i=1}^N$,
\begin{align*}
& \max_{|I|=j_{s+1}-1} \left(\sum_{i \in I} (\pi_{s+1} v)_i^2\right)^{1/2}
\\
\leq & \left(\sum_{i =j_s}^{j_{s+1}-1} ((\pi_{s+1}v)_i^*)^2\right)^{1/2}+\left(\sum_{i < j_s} ((\Delta_{s}v)^*_i)^2\right)^{1/2}
+ \max_{|I|=j_{s}-1} \left(\sum_{i \in I} (\pi_{s}v)_i^2\right)^{1/2}
\\
\leq & \nu_{s+1}d(V) + \|\Delta_{s} v\|_{[s]} +
\max_{|I|=j_{s}-1} \left(\sum_{i \in I} (\pi_{s}v)_i^2\right)^{1/2}.
\end{align*}
Thus, for every $v \in V$,
\begin{equation*}
\left(\sum_{i=1}^N (\pi_{s_1} v)_i^2\right)^{1/2}
\leq  d(V)\left(\sqrt{N}+\sum_{s=s_0}^{s_1-1} \nu_{s+1}\right) +  \sum_{s=s_0}^{s_1-1} \|\Delta_{s} v\|_{[s]} + \|\pi_{s_0}v\|_{[s_0]}.
\end{equation*}
Next, recall that $v =\sum_{s \geq s_1} \Delta_s v + \pi_{s_1}v$ and that $j_s=N+1$ for $s \geq s_1$. Therefore, by Assumption \ref{ass:quad-strucutre} for $p=1$,
$$
\left(\sum_{i=1}^N (\Delta_s v)_i^2 \right)^{1/2} \leq \|\Delta_s v\|_{[s]}.
$$
Combining the two estimates, one has
\begin{align*}
\left(\sum_{i=1}^N v_i^2\right)^{1/2} \leq & \sum_{s \geq s_1} \|\Delta_s v\|_{\ell_2^N} + \|\pi_{s_1} v\|_{\ell_2^N}
\\
\leq & \sum_{s \geq s_0} \|\Delta_s v\|_{[s]} + \|\pi_{s_0}v\|_{[s_0]} + d(V) \left(\sqrt{N}+\sum_{s=s_0}^{s_1-1} \nu_{s+1}\right)
\\
\leq & \Lambda(V)+d(V)\left(\sqrt{N}+\sum_{s=s_0}^{s_1-1} \nu_{s+1}\right).
\end{align*}
\endproof

With Lemma \ref{lemma:ell-2-N-radius} in place, let us turn our attention to the product Bernoulli process. Let $V,W \subset \R^N$ that satisfy Assumption \ref{ass:quad-strucutre} for $p=1$ and $p=2$ (and with the same admissible sequence for both values of $p$). Assume further that for the integers $s_0$ and $s_1$ as above,
$$
\sum_{s = s_0}^{s_1-1} \nu_{s+1} \leq c\sqrt{N}
$$
for a suitable absolute constant $c$.
\begin{Theorem} \label{thm:Bernoulli-prod}
There exist absolute constants $c_1$ and $c_2$ for which the following holds. If $V$, $W$, $j_s$ and $s_0$ are as above, then for every $t>4$, with probability at least $1-2\exp(-c_1t^2 2^{s_0})$,
\begin{align*}
\sup_{v \in V, \ w \in W} \left|\sum_{i=1}^N \eps_i v_i w_i \right| \leq & c_2\left(\Lambda(V)\Lambda(W) +\sqrt{N}\left(d(W)\Lambda(V)+d(V)\Lambda(W)\right)\right)
\\
+& c_2t\sqrt{N}\left(d(W)\Theta(V)+d(V)\Theta(W)\right).
\end{align*}
\end{Theorem}

\proof Observe that coordinate-wise,
\begin{align*}
\pi_{s+1} v \cdot \pi_{s+1} w = & (\pi_{s+1} v - \pi_s v) \cdot \pi_{s+1} w + \pi_s v \cdot (\pi_{s+1} w -\pi_s w) + \pi_s v \cdot \pi_s w
\\
= & \Delta_{s} v \cdot \pi_{s+1} w + \pi_s v \cdot \Delta_s w + \pi_s v \cdot \pi_s v.
\end{align*}
Therefore,
\begin{align*}
& \sup_{v \in V, \ w \in W} \left|\sum_{i=1}^N \eps_i v_i w_i \right| \leq  \sup_{v \in V, \ w \in W} \sum_{s \geq s_0}  \left| \sum_{i=1}^N \eps_i \left(\Delta_{s} v \cdot \pi_{s+1} w + \pi_s v \cdot \Delta_s w\right)_i \right|
\\
+ & \sup_{v \in V, \ w \in W} \left| \sum_{i=1}^N \eps_i (\pi_{s_0} v \cdot \pi_{s_0} w)_i \right|.
\end{align*}
By H\"{o}ffding's inequality, for a fixed $I \subset \{1,...,N\}$ and $z \in \R^N$, one has that with probability at least $1-2\exp(-t^2/2)$,
$$
\left|\sum_{i=1}^N \eps_i z_i \right| \leq \sum_{i \in I} |z_i| + t\left(\sum_{i \in I^c} z_i^2\right)^{1/2}.
$$
Note that by Assumption \ref{ass:quad-strucutre}, Lemma \ref{lemma:ell-2-N-radius} and since $j_{s-1} \leq j_s$,
$$
\left(\sum_{i < j_s} ((\Delta_s v)_i^*)^2\right)^{1/2} \leq \|\Delta_s v\|_{[s]}, \ \ \ \left(\sum_{i \geq j_s} ((\Delta_s v)_i^*)^4\right)^{1/4} \leq N^{1/4} \|\Delta_s v\|,
$$
$$
\left(\sum_{i=1}^N (\pi_{s+1}w)_i^2 \right)^{1/2} \leq \Lambda(W)+c\sqrt{N}d(W)
$$
and
$$
\left(\sum_{i \geq j_s} ((\pi_{s+1}w)_i^*)^4\right)^{1/4} \leq N^{1/4} d(W).
$$
Let $I$ be the union of the $j_s-1$ largest coordinates of $(|\Delta_s v|_i)_{i=1}^N$ and the $j_s-1$ largest coordinates of $(|\pi_{s+1}w|_i)_{i=1}^N$. Since $|\{\Delta_s v : v \in V\}| \leq 2^{2^{s+2}}$ and $|\{\pi_{s+1} w : w \in W\}| \leq 2^{2^{s+1}}$, it follows that for $t \geq 4$, with probability at least $1-2\exp(-ct^2 2^s)$, for every $v \in V$ and $w \in W$,
\begin{align*}
& \left|\sum_{i=1}^N \eps_i (\Delta_s v)_i (\pi_{s+1} w)_i \right| \leq  \sum_{i \in I} |(\Delta_s v)_i (\pi_{s+1} w)_i| + t2^{s/2} \left(\sum_{i \in I^c} (\Delta_s v)_i^2 (\pi_{s+1} w)_i^2 \right)^{1/2}
\\
\leq & \left(\sum_{i \in I} (\Delta_s v)_i^2\right)^{1/2} \left(\sum_{i \in I} (\pi_{s+1} w)_i^2\right)^{1/2} + t 2^{s/2} \left(\sum_{i \in I^c} (\Delta_s v)_i^4\right)^{1/4} \left(\sum_{i \in I^c} (\pi_{s+1} w)_i^4\right)^{1/4}
\\
\leq & \sqrt{2}\left(\sum_{i < j_s} ((\Delta_s v)_i^*)^2 \right)^{1/2}\left(\sum_{i=1}^{N}(\pi_{s+1} w)_i^2 \right)^{1/2} + t 2^{s/2} \left(\sum_{i \geq j_s} ((\Delta_s v)_i^*)^4\right)^{1/4} \left(\sum_{i \geq j_s}^N ((\pi_{s+1} w)_i^*)^4\right)^{1/4}
\\
\lesssim & \|\Delta_s v\|_{[s]} (\Lambda(W)+\sqrt{N}d(W)) + t2^{s/2} N^{1/4}\|\Delta_s v\| \cdot N^{1/4} d(W).
\end{align*}
Repeating the argument for $((\pi_s v)_i \cdot (\Delta_s w)_i)_{i=1}^N$ and summing over $s  \geq s_0$, one has that with probability at least $1-2\exp(-c_1t^2 2^{s_0})$, for every $v \in V$ and $w \in W$
\begin{align*}
& \sum_{s \geq s_0} \left| \sum_{i=1}^N \eps_i \left(\Delta_{s} v \cdot \pi_{s+1} w + \pi_s v \cdot \Delta_s w\right)_i \right|
\\
\lesssim &  \Lambda(V)\Lambda(W) + \sqrt{N}\left(d(W)\Lambda(V)+d(V)\Lambda(W)\right)
\\
+ &
t\sqrt{N}\left(d(W)\Theta(V)+d(V)\Theta(W)\right).
\end{align*}

Applying the same type of argument to
$$
(*)=\left|\sum_{i=1}^N \eps_i (\pi_{s_0}v)_i \cdot (\pi_{s_0}w)_i\right|,
$$
it is evident that with probability at least $1-2\exp(-c_2t^22^s)$,
\begin{align*}
(*) \lesssim & \|\pi_{s_0}v\|_{[s_0]}(\Lambda(W)+\sqrt{N}d(W))
\\
+ & t 2^{s_0/2} \cdot N^{1/4}\|\pi_{s_0}v\| \cdot N^{1/4}\|\pi_{s_0}w\|
\\
\lesssim & \Lambda(V) (\Lambda(W)+\sqrt{N}d(W)) + t\sqrt{N}d(V)\Theta(W),
\end{align*}
which concludes the proof.
\endproof

Corollary \ref{cor:multi-in-section} and Theorem \ref{thm:Bernoulli-prod} are the first component in the proofs of Theorem \ref{thm:multiplier-intro} and Theorem \ref{thm:quadratic-process-upper-intro}, respectively. For the other component, one has to show that typical coordinate projections of the indexing classes are well-behaved in the sense of Assumption \ref{ass:decomp-of-a-set} or of Assumption \ref{ass:quad-strucutre}. To that end, one has to identify the norms $\| \ \|_{[s]}$ and $\| \ \|$, the sequences $(j_s)_{s \geq s_0}$ and $(\nu_s)_{s \geq s_0}$ and estimate the resulting complexity terms, $\Lambda(P_\sigma F)$, $\Theta(P_\sigma F)$ and $d(P_\sigma F)$. The main step towards that goal is presented in the next section.

\section{Structural results - preliminary estimates} \label{sec:pre}
Let $Z$ be a random variable. Our primary goal is to study the monotone nonincreasing rearrangement of $N$ independent copies of $Z$. We will show that the norms 
$$
\|Z\|_{(p)} = \sup_{1 \leq q \leq p} \frac{\|Z\|_{L_q}}{\sqrt{q}}
$$
play a key role in the desired estimate.
\begin{Theorem} \label{thm:decom-for-iid}
There exist absolute constants $c_0,c_1,c_2$ for which the following holds. Let $1 \leq r<q$ and set $0<\beta<(q/r)-1$. Consider $Z \in L_q$, and let $Z_1,...,Z_N$ be independent copies of $Z$. Put $1 \leq p \leq N$ and set
$$
j_0 = \min\left\{\left\lceil \frac{c_0p}{((q/r)-1) \log(4+eN/p)}\right\rceil,N+1\right\}.
$$
The random vector $(Z_i)_{i=1}^N$ can be decomposed to a sum $U+V$, for random vectors $U,V \in \R^N$ that have disjoint supports and satisfy:

If $j_0 \geq 2$, then
\begin{description}
\item{$\bullet$} $|{\rm supp}(U)| = j_0-1$.
\item{$\bullet$} For every $t>2$ with probability at least $1-t^{-2p}\exp(-p)$,
$$
\|U\|_{\ell_2^N} \leq c_1t \sqrt{2p}\|Z\|_{(2p)}.
$$
\item{$\bullet$} For every $t>2$, with probability at least $1-t^{-j_0q}\exp(-p)$,
$$
\|V\|_{\ell_r^N} \leq c_1\left(\frac{q}{q-r}\right)^{1/r} t \|Z\|_{L_q} N^{1/r}.
$$
\end{description}
And, if $j_0=1$ then
\begin{description}
\item{$\bullet$} $U=0$.
\item{$\bullet$} For $t>2$ with probability at least $1-c_2t^{-q}N^{-\beta}$,
$$
\|V\|_{\ell_r^N} \leq c_1 \left(\frac{q}{q-(\beta+1)r}\right)t \|Z\|_{L_q} N^{1/r}.
$$
\end{description}
\end{Theorem}

We will show that the vector $U$ consists of the $j_0-1$ largest coordinates of $(|Z_i|)_{i=1}^N$ and that $V=(Z_i)_{i=1}^N-U$. The key is to determine correctly the `cut-off' point in that decomposition -- which happens to be $j_0$.

\vskip0.5cm

As the formulation of Theorem \ref{thm:decom-for-iid} indicates, the treatment depends on whether or not the decomposition is trivial (i.e., if $U=0$), and on the required probability estimate.

\vskip0.5cm

We begin with the smaller coordinates of $(|Z_i|)_{i=1}^N$, which will be used to define the vector $V$ in the decomposition of $(Z_i)_{i=1}^N$.

\begin{Lemma} \label{lemma:single-small-coordinates}
There exist absolute constants $c_0,c_1$ and $c_2$ for which the following holds. Let $1 \leq r < q$, set $Z \in L_q$ and put $Z_1,...,Z_N$ to be independent copies of $Z$. Fix $1 \leq p \leq N$,  let 
$$
j_0 = \left\lceil \frac{c_0p}{((q/r)-1) \log(4+eN/p)}\right\rceil
$$ 
and $t>2$, as above.

If $j_0>1$, then with probability at least $1-2t^{-j_0q}\exp(-p)$,
$$
\left(\sum_{j=j_0}^{N} (Z_i^*)^r \right)^{1/r} \leq c_1\left(\frac{ q}{q-r}\right)^{1/r} t N^{1/r} \|Z\|_{L_q}.
$$
And, if $j_0=1$ and $0<\beta<(q/r)-1$ then with probability at least $1-c_2t^{-q}N^{-\beta}$,
$$
\left(\sum_{j=1}^{N} |Z_i|^r \right)^{1/r} \leq c_1 \left(\frac{q}{q-(\beta+1)r}\right)t \|Z\|_{L_q} N^{1/r}.
$$
\end{Lemma}

\proof
Fix $\rho \geq 1$ to be named later and observe that by a binomial estimate and Markov's inequality, for every $u>0$,
$$
Pr(Z_j^*>u) \leq \binom{N}{j}Pr^j(|Z|>u) \leq \left(\frac{eN \|Z\|_{L_q}^q}{ju^q}\right)^j.
$$
Therefore, if $u=t \|Z\|_{L_q} (eN/j)^{\rho/q}$ then for any $j \leq N$,
\begin{align} \label{eq:binomial-est-inproof}
& Pr\left( \exists k \geq j, \ Z_k^* \geq t \|Z\|_{L_q} \left(\frac{eN}{k}\right)^{\rho/q}\right) \leq \sum_{k \geq j} \frac{1}{t^{kq}} \left(\frac{eN}{k}\right)^{-k(\rho-1)}
\nonumber
\\
\leq & \frac{2}{t^{j q}} \cdot \left(\frac{eN}{j}\right)^{-j(\rho-1)}.
\end{align}

Note that $\left({eN}/{j}\right)^{-j(\rho-1)} \leq \exp(-p)$ when
$$
j \geq j_0 = \left\lceil \frac{c_2 p}{ (\rho-1) \log(4+e(\rho-1)N/p)}\right\rceil,
$$
and on that event \eqref{eq:binomial-est-inproof},
\begin{equation*}
\sum_{j=j_0}^{N} (Z_j^*)^r \leq t^{r} \|Z\|_{L_q}^r \sum_{j = j_0}^{N} \left(\frac{e N}{j}\right)^{\rho r/q}.
\end{equation*}
Set $\alpha=\rho r/q$ and observe that if $\alpha<1$, that is, if $\rho < q/r$, then
$$
\sum_{j =j_0}^{N} (Z_j^*)^r \leq \frac{c_3}{1-\alpha} t^{r} \|Z\|_{L_q}^r N
$$
for an absolute constant $c_3$. The first part of the claim now follows by setting $\rho=1+(q/r-1)/2$, implying that $1-\alpha=(q-r)/2q$.

The second part follows an identical path to the first one, by setting $0<\beta<(q/r)-1$, $\rho=\beta+1$ and $\alpha=(\beta+1)r/q<1$.
\endproof

\begin{Remark} \label{rem:small-coordinates}
Let $q>4$ and $r=2$. Set $\beta=1/2 < (q/r)-1$ and thus $\alpha<3/4$. As noted in the proof of Lemma \ref{lemma:single-small-coordinates}, with probability at least $1-\exp(-p)$, for every $i \geq j_0$
$$
Z_i^* \lesssim \|Z\|_{L_q} \left(\frac{eN}{i}\right)^{\alpha/2}.
$$
Set $a_i=(eN/i)^{\alpha/2}$ and consider an increasing sequence $j_0 < j_1 < ... < j_\ell = N$ for which
$$
\sum_{k=0}^{\ell} j_k^{1/2-\alpha/2} \leq c_1N^{(1-\alpha)/2}.
$$
It is straightforward to verify that
\begin{equation} \label{eq:sum-of-monotone}
\sum_{k=0}^{\ell-1} \left(\sum_{i=j_k}^{j_{k+1}} a_i^2\right)^{1/2} \leq c_2\sqrt{N},
\end{equation}
and $c_2$ depends only on $c_1$.

This simple fact will play a role in verifying the third part of Assumption \ref{ass:quad-strucutre} for a typical coordinate projection.
\end{Remark}

Next, let us show how the norms $\| \ \|_{(p)}$ may be used to upper bound the larger coordinates in a monotone rearrangement of $Z_1,...,Z_N$.

\begin{Lemma} \label{lemma:large-coordinates}
Let $Z_1,...,Z_N$ be independent copies of a random variable $Z$, set $p \geq \log N$ and put $1 \leq m \leq N/2e$ for which $\binom{N}{m} \leq \exp(p)$. Then, for every $t>1$, with probability at least $1-t^{-2p}\exp(-p)$, one has
$$
\left(\sum_{i \leq m} (Z_i^*)^2\right)^{1/2} \lesssim  t\sqrt{p}\|Z\|_{(2p)}.
$$
\end{Lemma}

The proof of Lemma \ref{lemma:large-coordinates} is based on the following fact, due to Lata{\l}a \cite{Lat-moments}.
\begin{Theorem} \label{thm:latala}
Let $W$ be a nonnegative random variable. If $W_1,...,W_m$ are independent copies of $W$, then
$$
\|\sum_{i=1}^m W_i \|_{L_r} \sim \sup\left\{\frac{r}{s}\left(\frac{m}{r}\right)^{1/s} \|W\|_{L_s} \ \ : \ \ \max\left\{1,\frac{r}{m}\right\} \leq s \leq r \right\}.
$$
\end{Theorem}

\noindent{\bf Proof of Lemma \ref{lemma:large-coordinates}.}
Since $\binom{N}{m} \geq (\frac{N}{m}-1)^m \geq \exp(m)$, it follows that $p \geq m$; thus, by Theorem \ref{thm:latala} for $W=Z^2$ and $r=p$,
$$
\|\sum_{i=1}^m Z_i^2\|_{L_{p}} \lesssim \max\left\{ \frac{p}{s} \|Z^2\|_{L_s}  \ : \ \frac{p}{m} \leq s \leq p \right\}.
$$
Clearly, for $s \leq p$,
$$
\|Z^2\|_{L_s} = \|Z\|_{L_{2s}}^2 \leq 2s \|Z\|_{(2p)}^2,
$$
and $\|\sum_{i=1}^m Z_i^2\|_{L_{p}} \lesssim 2p\|Z\|_{(2p)}^2$. Therefore,
\begin{align*}
& Pr\left( \exists I \subset \{1,...,N\}, \ |I|=m \ : \left(\sum_{i \in I} Z_i^2\right)^{1/2} \geq eu \right)
\\
\leq & \binom{N}{m} Pr\left(\left(\sum_{i=1}^m Z_i^2\right) \geq e^2u^2 \right)
\\
\leq & \binom{N}{m} \frac{\|\sum_{i=1}^m Z_i^2\|_{L_{p}}^{p}}{(eu)^{2p}} \leq \binom{N}{m} \left(\frac{cp \|Z\|_{(2p)}^2}{e^2u^2}\right)^p \leq t^{-2p}\exp(-p)
\end{align*}
for the choice of $u \sim t\sqrt{2p}\|Z\|_{(2p)}$ and since $\binom{N}{m} \leq \exp(p)$.
\endproof

\vskip0.5cm

The proof of Theorem \ref{thm:decom-for-iid} is simply the combination of Lemma \ref{lemma:single-small-coordinates} and Lemma \ref{lemma:large-coordinates}.
\endproof

\vskip0.5cm

Let us turn to the main application of Theorem \ref{thm:decom-for-iid}.

\vskip0.5cm

From here on, fix $u \geq 4$ and for every $s \geq 0$ and $r<q$ set
\begin{equation} \label{eq:j-s-def}
j_s(r,q)=\min\left\{\left\lceil\frac{c_0u^22^s}{((q/r)-1)\log(4+eN/u^22^s)}\right\rceil,N+1\right\}
\end{equation}
for a suitable absolute constant $c_0$ as in Theorem \ref{thm:decom-for-iid}. Consider a finite class of functions $H$, whose cardinality is at most $2^{2^{s+3}}$.

Writing $j_s$ instead of $j_s(r,q)$, let us examine three different cases: $j_s = N+1$, $2 \leq j_s \leq N$ and $j_s=1$. Motivated by the requirements of the chaining arguments outlined earlier, in all three cases one would like to obtain uniform control over all the functions in $H$; thus, the probability estimate with which one must control the decomposition from Theorem \ref{thm:decom-for-iid} for each individual function should be at least $1-\exp(-2^{s+3})$.

\vskip0.5cm
\begin{description}
\item{$\bullet$} When $j_s=N+1$, the decomposition is trivial, in the sense that for each function $h$, $U=(h(X_i))_{i=1}^N$ and $V=0$. Hence, setting $2p=u^22^s$, it follows that with probability at least $1-\exp(-u^22^s/2)$, for every $h \in H$, $(h(X_i))_{i=1}^N=U$, and
$$
\|U\|_{\ell_2^N} \lesssim  u2^{s/2} \|h\|_{(u^22^s)}.
$$
\item{$\bullet$} When $1 < j_s \leq N$, and setting $2p=u^22^s$ once again, it follows that with probability at least $1-2\exp(-u^22^s/2)$, for every $h \in H$, $(h(X_i))_{i=1}^N = U+V$, where
$$
\|U\|_{\ell_2^N} \lesssim  u2^{s/2} \|h\|_{(u^22^s)}, \ \ {\rm and} \ \ \|V\|_{\ell_r^N} \leq c(q,r)\|h\|_{L_q} N^{1/r}.
$$
\item{$\bullet$} When $j_s=1$, $(h(X_i))_{i=1}^N = V$ and $U=0$. Moreover, because $j_s=1$,
$$
c_0u^22^s \leq ((q/r)-1)\log(eN/u^22^s);
$$
hence,
$$
|H| \leq 2^{2^{s+3}} \leq (c_1 N)^{c_2/u^2},
$$
for constants $c_1$ and $c_2$ that depend only on $r$ and $q$.

Let $0<\beta \leq (q/r)-1$ and note that by Theorem \ref{thm:decom-for-iid}, with probability at least $1-c_3N^{-\beta}$,
\begin{equation} \label{eq:triangle}
\|V\|_{\ell_r^N} \leq c_4 \left(\frac{q}{q-(\beta+1)r}\right) \|h\|_{L_q} N^{1/r}.
\end{equation}
Therefore, \eqref{eq:triangle} holds with probability of $1-2\exp(-\theta u^2 2^s)$ for every $h \in H$ if
\begin{equation} \label{eq:prob-est-finite-class}
c_3N^{-\beta} 2^{2^{s+3}} \leq \exp(-\theta u^2 2^s),
\end{equation}
which is the case when $u \geq c_4(q,r)/\sqrt{\beta}$ and $\theta \leq c_5(q,r)\beta$.
\end{description}
Combining these observations yields the following outcome:
\begin{Corollary} \label{cor:decomposition-finite-class}
There exist absolute constants $c_0,c_1$ and for every $1 \leq r <q$ there exist constants $c_2$ and $c_3$ that depend only on $q$ and $r$ for which the following holds. Set
$$
j_s=\min\left\{\left\lceil\frac{c_0u^22^s}{((q/r)-1)\log(4+eN/u^22^s)}\right\rceil,N+1\right\}
$$
for $0<\beta \leq (q/r)-1$ and $u \geq c_2/\sqrt{\beta}$. If $H \subset L_q$ is of cardinality at most $2^{2^{s+2}}$, then with probability at least $1-2\exp(-c_3\beta u^2 2^s)$, for every $h \in H$, $(h(X_i))_{i=1}^N=U_h+V_h$; the support of each $U_h$ is the set of the largest $j_s-1$ coordinates of $(|h(X_i)|)_{i=1}^N$ while $V_h$ is supported on its complement;
$$
\|U_h\|_{\ell_2^N} \leq c_1u2^{s/2}\|h\|_{(u^22^s)}; \ \ {\rm and} \ \ \|V_h\|_{\ell_r^N} \leq c_1\left(\frac{q}{q-(\beta+1)r}\right)\|h\|_{L_q} N^{1/r}.
$$
\end{Corollary}

\section{Proofs of the main results} \label{sec:main}
Let us turn to the implications of Corollary \ref{cor:decomposition-finite-class} in the contexts of Assumption \ref{ass:decomp-of-a-set} and Assumption \ref{ass:quad-strucutre}.

\vskip0.5cm

Let $F$ be a class of functions and set $(F_s)_{s \geq 0}$ to be an admissible sequence of $F$. Fix $s_0$ to be named later and set $s \geq s_0$. Clearly,
$$
|F_s|, \ |\left\{\Delta_s f : f \in F\right\}| \leq 2^{2^{s+2}},
$$
and thus Corollary \ref{cor:decomposition-finite-class} holds for both $F_s$ and $\{\Delta_s f : f \in F\}$.

For every $s \geq s_0$, denote by $\left[(\Delta_s f)(X_i)\right]^*$ the $i$-th largest coordinate of the vector $\left(|\Delta_s f|(X_i)\right)_{i=1}^N$, and put $\left[(\pi_{s} f)(X_i)\right]^*$ to be the $i$-th largest coordinate of the vector $\left(|\pi_{s} f|(X_i)\right)_{i=1}^N$. Finally, recall that
$$
\tilde{\Lambda}_{s_0,u}(F) \leq \Lambda_{s_0,u}(F)+2^{s_0/2}\sup_{f \in F} \|f\|_{(u^22^{s_0})},
$$

\subsection{The multiplier process} \label{sec:multi}
Let us see how Corollary \ref{cor:decomposition-finite-class} may be used to prove Theorem \ref{thm:multiplier-intro}.

Let $q>2$, $r=\min\{1/2+q/4,2\}$ and $r_1=2r^\prime$, where $r^\prime$ is the  conjugate index of $r$. Put $q_1=2r_1$ and set $\beta=1/2 < (q_1/r_1)-1$. Also, let
$$
j_s=j_s(r_1,q_1)=\min \left\{\left\lceil \frac{c_0u^2 2^s}{\log(4+eN/2^s)}\right\rceil,N+1\right\}.
$$

Below is the summary of the outcome of Corollary \ref{cor:decomposition-finite-class} when applied to the classes $\{\Delta_s f : f \in F\}$ and $F_s$ for $s \geq s_0$, $q_1$ and $r_1$ as above, and for $u > \max\{\sqrt{q_1},4\}$.
\begin{Corollary} \label{cor:decomp-for-multipliers}
There are constants $c_1$ and $c_2$ that depend only on $q$ and an event of probability at least $1-2\exp(-c_1u^2 2^{s_0})$ on which the following holds. For every $f \in F$ and $s \geq s_0$,
$$
\left( \sum_{i < j_s} \left[(\Delta_s f)^2(X_i)\right]^* \right)^{1/2} \leq c_2u 2^{s/2} \|\Delta_s f\|_{(u^22^s)},
$$
$$
\left( \sum_{i \geq j_{s}} \left[(\Delta_s f)^{2r^\prime}(X_i)\right]^* \right)^{1/2r^\prime} \leq c_2 N^{1/2r^\prime} \|\Delta_s f\|_{L_{q_1}},
$$
$$
\left( \sum_{i < j_s} \left[(\pi_{s} f)^2(X_i)\right]^* \right)^{1/2} \leq c_2 u 2^{s/2} \|\pi_{s} f\|_{(u^22^{s})},
$$
and
$$
 \left( \sum_{i \geq j_{s}} \left[(\pi_{s} f)^{2r^\prime}(X_i)\right]^* \right)^{1/2r^\prime} \leq c_2 N^{1/2r^\prime} \|\pi_{s} f\|_{L_{q_1}}.
$$
\end{Corollary}

\begin{Remark} \label{rem-multiplier-large-q}
Observe that if $q \geq 4$ then $r=2$, and thus $2r=2r^\prime=4$ and $q_1=8$. Therefore, all the constants in Corollary \ref{cor:decomp-for-multipliers} are absolute constants and one may take any $u \geq 8$.
\end{Remark}

Corollary \ref{cor:decomp-for-multipliers} implies that with probability at least $1-2\exp(-c_1u^2 2^{s_0})$, the coordinate projection $P_\sigma F \subset \R^N$ satisfies Assumption \ref{ass:decomp-of-a-set} -- with the natural identification of elements in $P_\sigma F$ with functions in $F$ via the coordinate projection map, and for an admissible sequence in $P_\sigma F$ endowed by one in $F$ -- for the choices of
\begin{description}
\item{$\bullet$} $\| \ \|_{[s]} = c_2u2^{s/2}\| \ \|_{(u^22^s)}$,
\item{$\bullet$} $\| \ \| = c_2\| \ \|_{L_{q_1}}$,
\item{$\bullet$} $p=r^\prime$
\end{description}
Hence, for an almost optimal admissible sequence in $F$,
\begin{align*}
\Lambda(P_\sigma F) \lesssim & u \left(\Lambda_{s_0,u}(F)+2^{s_0/2}\sup_{f \in F}\|f\|_{(u^22^s)}\right), \ \ {\rm and}
\\
\Theta(P_\sigma F) \lesssim & \gamma_{s_0,2}(F,L_{q_1})+2^{s_0/2}\sup_{f \in F} \|f\|_{L_{q_1}}.
\end{align*}

Recall that $u^2 \geq q_1$, and thus, for every $s \geq 0$,
$$
\frac{\|f\|_{L_{q_1}}}{\sqrt{q_1}} \leq \sup_{1 \leq q \leq u^22^s} \frac{\|f\|_{L_q}}{\sqrt{q}}=
\|f\|_{(u^22^s)}.
$$
Therefore,
\begin{align*}
\gamma_{s_0,2}(F,L_{q_1})+2^{s_0/2}\sup_{f \in F} \|f\|_{L_{q_1}} \lesssim & \sqrt{q_1} \left(\Lambda_{s_0,2}(F)+2^{s_0/2}\sup_{f \in F}\|f\|_{(u^22^s)}\right)
\\
= & \sqrt{q_1} \tilde{\Lambda}_{s_0,u}(F).
\end{align*}
Combining these observations with Corollary \ref{cor:multi-in-section}, it follows that for every $(X_i)_{i=1}^N$ for which Corollary \ref{cor:decomp-for-multipliers} holds (i.e., with probability at least $1-2\exp(-c_1u^22^{s_0})$ with respect to $\mu^N$), and for every $(\xi_i)_{i=1}^N$, one has that with $(\eps_i)_{i=1}^N$ probability at least $1-2\exp(-c_1t^22^{s_0})$, for every $f  \in F$
\begin{equation*}
\left|\sum_{i=1}^N \eps_i \xi_i f(X_i) \right| \lesssim
\tilde{\Lambda}_{s_0,u}(F) \cdot \left(u\|(\xi_i)_{i=1}^N\|_{\ell_2^N} +  t\sqrt{q_1}N^{1/2r^\prime} \left(\sum_{i \geq j_{s_0}} (\xi_i^*)^{2r}\right)^{1/2r}\right).
\end{equation*}
Hence, if
$$
\|(\xi_i)_{i=1}^N\|_{\ell_2^N} \leq A\sqrt{N}, \ \  \ \ \left(\sum_{i \geq j_{s_0}} (\xi_i^*)^{2r}\right)^{1/2r} \leq B N^{1/2r},
$$
and setting $t=u$, then
$$
\left|\sum_{i=1}^N \eps_i \xi_i f(X_i) \right| \lesssim u\tilde{\Lambda}_{s_0,u}(F) \cdot \sqrt{N}(A+\sqrt{q_1}B).
$$
Thus, to conclude the proof of Theorem \ref{thm:multiplier-intro}, one has identify $A$ and $B$ for which, with high probability,
$$
\|(\xi_i)_{i=1}^N\|_{\ell_{2}^N} \leq A \sqrt{N} \ \ {\rm and} \ \ \left(\sum_{i \geq j_{s_0}} (\xi_i^*)^{2r}\right)^{1/2r} \leq B N^{1/2r},
$$
and then apply the symmetrization argument of Theorem \ref{thm:GZ}.

\vskip0.5cm

By Lemma \ref{lemma:single-small-coordinates} and since $2r<1+q/2<q$, one has that with probability at least $1-2\exp(-c_0u^22^{s_0})$,
$$
\left(\sum_{i \geq j_{s_0}} (\xi_i^*)^{2r}\right)^{1/2r} \leq c(q) \|\xi\|_{L_q} N^{1/2r}.
$$
Therefore, one may take $B \sim_q \|\xi\|_{L_q}$ and all that remains is to deal with $\|(\xi_i)_{i=1}^N\|_{\ell_2^N}$. To that end, the following is a minor modification of Lemma \ref{lemma:single-small-coordinates}.

\begin{Lemma} \label{lemma:random-vector}
Let $q>2$ and assume that $\xi \in L_q$. If $\xi_1,...,\xi_N$ are independent copies of $\xi$ and $z=(\xi_i)_{i=1}^N$, then for every $w>1$, with probability at least $1-c_0w^{-q}N^{-((q/2)-1)}\log^{q} N$,
$$
\|z\|_{\ell_{2}^N} \leq c_1w \|\xi\|_{L_q}\sqrt{N},
$$
where $c_0$ and $c_1$ depend only on $q$.
\end{Lemma}

\proof
Let $\eta=(q/2)-1$, fix $1 \leq k \leq N/2$ and set $v>0$ to be named later. A binomial estimate implies that
\begin{align*}
& Pr(\xi_k^* \geq v(eN/k)^{(1+\eta)/q} \|\xi\|_{L_q}) \leq \binom{N}{k} Pr^k (|\xi| \geq v(eN/k)^{(1+\eta)/q}\|\xi\|_{L_q})
\\
\leq & \left(\frac{eN}{k}\right)^k \left(\frac{k}{eN}\right)^{(1+\eta)k} \cdot v^{-kq} =\left(\frac{e N}{k}\right)^{-\eta k} v^{-k q}.
\end{align*}
Hence, for $v=u (eN/k)^{1/2 - (1+\eta)/q}$,
$$
\xi_k^* \lesssim u\|\xi\|_{L_q}\sqrt{N/k}
$$
with probability at least
$$
1-u^{-kq}\left(\frac{ek}{N}\right)^{k((q/2)-1)}.
$$
For every $1 \leq k \leq N$, set $u_k=w/\log(eN/k)$ and observe that
$$
\sum_{k=1}^N \xi_k^2 \lesssim w^2N \|\xi\|_{L_q}^{2} \sum_{k=1}^N \frac{1}{k\log^2(eN/k)} \leq c_1 w^2N\|\xi\|_{L_q}^{2}.
$$
The claim follows by summing the probability estimates.
\endproof

Lemma \ref{lemma:random-vector} implies that one may select $A \sim_q w \|\xi\|_{L_q}$ and with probability at least $1- c_0(q)w^{-q}N^{-((q/2)-1)}\log^{q} N-2\exp(-c_1u^22^{s_0})$,
\begin{equation} \label{eq:multi-proof-end}
\left|\sum_{i=1}^N \eps_i \xi_i f(X_i) \right| \lesssim_q uw\sqrt{N} \|\xi\|_{L_q} \cdot \tilde{\Lambda}_{s_0,u}(F).
\end{equation}
\endproof

Let us turn to a version of Theorem \ref{thm:multiplier-intro} when $\xi \in L_{\psi_2}$.

\begin{Theorem} \label{thm:multiplier-psi-2}
There exist absolute constants $c_1$ and $c_2$ for which the following holds. If $\xi \in L_{\psi_2}$ then for every $u,w \geq 8$, with probability at least $1-2\exp(-c_1u^22^{s_0})-2\exp(-c_1N w^2)$,
$$
\sup_{f \in F} \left|\sum_{i=1}^N \left(\xi_i f(X_i) - \E \xi f\right)\right| \leq c_2 uw \sqrt{N} \|\xi\|_{\psi_2}\tilde{\Lambda}_{s_0,u}(F).
$$
\end{Theorem}
The proof follows a similar path to the proof of Theorem \ref{thm:multiplier-intro} with a minor modification in the last step -- the bounds on $(\xi_i)_{i=1}^N$.

By Bernstein's inequality, with probability at least $1-2\exp(-c_0N \min\{w^2,w^4\})$,
$$
\left(\frac{1}{N}\sum_{i=1}^N \xi_i^2\right)^{1/2} \leq c_1(w+1)\|\xi\|_{\psi_2}.
$$
Therefore, if $u \geq 8$ and $w \geq 1$, then with probability at least
$$
1-2\exp(-c_0N w^2)-2\exp(-c_1u^22^{s_0}),
$$
$$
\left(\sum_{i=1}^N \xi_i^2\right)^{1/2} \leq c_2w\|\xi\|_{\psi_2}\sqrt{N} \ \ {\rm and} \ \ \left(\sum_{i \geq j_{s_0}} (\xi_i^*)^4\right)^{1/4} \leq c_2u\|\xi\|_{\psi_2} N^{1/4}
$$
for absolute constants $c_0,c_1$ and $c_2$.

The rest of the proof is unchanged, for the choices of $r=r^\prime=2$ and $q_1=4r^\prime=8$, as noted in Remark \ref{rem-multiplier-large-q}.
\endproof

\subsection{The quadratic process} \label{sec:quad}
Following the same path as in the previous section, and thanks to Theorem \ref{thm:Bernoulli-prod}, one has to show that typical coordinate projections $P_\sigma F$ and $P_\sigma H$ satisfy Assumption \ref{ass:quad-strucutre} for $p=1$ and $p=2$.

Fix $q>4$ and let $j_s=j_s(2,q)$ be as in \eqref{eq:j-s-def}. Set
$$
\nu_s = \left(\sum_{i=j_s}^{j_{s+1}-1} \left(\frac{eN}{i}\right)^\alpha \right)^{1/2}
$$
for $\alpha<3/4$ as in Remark \ref{rem:small-coordinates}. It is straightforward to verify that if $s_1$ is the smallest integer for which $j_s = N+1$ then
$$
\sum_{s=s_0}^{s_1-1} \nu_s \leq c\sqrt{N}
$$
for an absolute constant $c$.

To handle the first and second parts of Assumption \ref{ass:quad-strucutre}, one may apply Corollary \ref{cor:decomposition-finite-class} to the classes
$\{\Delta_s f : f \in F\}$, $\{\Delta_s h : h \in H\}$, $F_s$ and $H_s$ for $s \geq s_0$.

\begin{Corollary} \label{cor:decomp-for-products}
There exists an absolute constant $c_1$, a constant $c_2$ that depends only on $q$ and an event of probability at least $1-2\exp(-c_1u^2 2^{s_0})$ on which the following holds. Consider $u \geq \sqrt{q}$ and set $j_s=j_s(2p,q)$ as in \eqref{eq:j-s-def}, for $p=1$ and $p=2$. For every $f \in F$ and every $s \geq s_0$,
$$
\left( \sum_{i < j_s} \left[(\Delta_s f)^2(X_i)\right]^* \right)^{1/2} \leq c_2u 2^{s/2} \|\Delta_s f\|_{(u^22^s)},
$$
$$
\left( \sum_{i \geq j_{s-1}} \left[(\Delta_s f)^{2p}(X_i)\right]^* \right)^{1/2p} \leq c_2 N^{1/2p} \|\Delta_s f\|_{L_{q}},
$$
$$
\left( \sum_{i < j_s} \left[(\pi_{s} f)^2(X_i)\right]^* \right)^{1/2} \leq c_2 u 2^{s/2} \|\pi_{s} f\|_{(u^22^{s})},
$$
and
$$
 \left( \sum_{i \geq j_{s-1}} \left[(\pi_{s} f)^{2p}(X_i)\right]^* \right)^{1/2p} \leq c_2 N^{1/2p} \|\pi_{s} f\|_{L_{q}}.
$$
\end{Corollary}

Therefore, with probability at least $1-2\exp(-c_1u^2 2^{s_0})$, the coordinate projection $P_\sigma F$ satisfies Assumption \ref{ass:quad-strucutre} for $p=1$ and $p=2$, with the choices of
\begin{description}
\item{$\bullet$} $\| \ \|_{[s]} = c_2u2^{s/2}\| \ \|_{(u^22^s)}$;
\item{$\bullet$} $\| \ \| = c_2\| \ \|_{L_{q}}$ (and, in particular, $d(P_\sigma F) \leq c_2\sup_{f \in F} \|f\|_{L_q}$);
\item{$\bullet$} $\nu_s = \left(\sum_{i=j_s}^{j_{s+1}-1} \left(\frac{eN}{i}\right)^\alpha \right)^{1/2}$, implying that $\sum_{s=s_0}^{s_1-1} \nu_s \leq c\sqrt{N}$.
\end{description}
Moreover, if $(F_s)_{s \geq s_0}$ is an almost optimal admissible sequence,
\begin{align*}
\Lambda(P_\sigma F) \lesssim & u \left(\Lambda_{s_0,u}(F)+2^{s_0/2}\sup_{f \in F}\|f\|_{(u^22^s)}\right) = u\tilde{\Lambda}_{s_0,u}(F),
\\
\Theta(P_\sigma F) \lesssim & \gamma_{s_0,2}(F,L_q)+2^{s_0/2}\sup_{f \in F} \|f\|_{L_{q}},
\end{align*}
and
$$
\gamma_{s_0,2}(F,L_q)+2^{s_0/2}\sup_{f \in F} \|f\|_{L_{q}} \lesssim \sqrt{q}\tilde{\Lambda}_{s_0,2}(F).
$$

This observation, together with Theorem \ref{thm:Bernoulli-prod} and the symmetrization argument of Theorem \ref{thm:GZ} completes the proof of Theorem \ref{thm:quadratic-process-upper-intro}.

\vskip0.5cm

\subsection{Unconditional log-concave ensembles}
We end this article with an example that shows yet again the advantage $\tilde{\Lambda}$ has over a $\psi_2$-based complexity term.

Let $X=(x_1,...,x_n)$ be a random vector, distributed according to an isotropic, unconditional, log-concave measure on $\R^n$. By that we mean that for every $x \in \R^n$, $\E \inr{x,X}^2 = \|x\|_{\ell_2^n}^2$; that $X$ has the same distribution as $(\eps_1x_1,...,\eps_nx_n)$ for every choice of signs $\eps_1,...,\eps_n$; and that $X$ has a log-concave density.

Applying the Bobkov-Nazarov Theorem \cite{BN}, $X$ is stochastically dominated by $Y=(y_1,...,y_n)$, a random vector whose coordinates are independent, standard, exponential random variables. Hence, by \cite{GluKwa}, for every $t \in \R^n$ and every $p \geq 2$,
$$
\|\inr{X,t}\|_{L_p} \leq \|\inr{Y,t}\|_{L_p} \sim p\|t\|_{\ell_\infty^n}+\sqrt{p}\|t\|_{\ell_2^n}.
$$
Therefore, $\|\inr{X,t}\|_{(u^22^s)} \lesssim u2^{s/2} \|t\|_{\ell_\infty^n} + \|t\|_{\ell_2^n}$ and for every $s \geq 0$,
$$
\Lambda_{s,u}(F_T) \leq c_1\left(u\gamma_1(T,\ell_\infty^n)+\gamma_2(T,\ell_2^n)\right).
$$
Set $s_0$ to be the largest integer for which $\gamma_1(T,\ell_\infty^n) \geq 2^s \sup_{t \in T} \|t\|_{\ell_\infty^n}$ and $\gamma_2(T,\ell_2^n) \geq 2^{s/2} \sup_{t \in T} \|t\|_{\ell_2^n}$. Hence,
\begin{align*}
2^{s_0/2}\sup_{t \in T} \|\inr{Y,t}\|_{(u^22^{s_0})} \lesssim & u2^{s_0} \sup_{t \in T}\|t\|_{\ell_\infty^n} + 2^{s_0/2}\sup_{t \in T}\|t\|_{\ell_2^n}
\\
\lesssim & u\gamma_1(T,\ell_\infty^n)+\gamma_2(T,\ell_2^n),
\end{align*}
and
$$
\tilde{\Lambda}_{s_0,u}(F_T) \lesssim u\gamma_1(T,\ell_\infty^n)+\gamma_2(T,\ell_2^n).
$$
As noted earlier, setting $E(T)=\E \sup_{t \in T } \inr{Y,t}$, then by a result due to Talagrand \cite{Tal-AM,Tal-book},
$$
\gamma_1(T,\ell_\infty^n)+\gamma_2(T,\ell_2^n) \sim \E \sup_{t \in T} \inr{Y,t},
$$
implying that for $u \geq 1$,
\begin{equation*}
\tilde{\Lambda}_{s_0,u}(F_T) \leq c_2u E (T).
\end{equation*}
Therefore, with probability at least $1-2\exp(-c_3u^22^{s_0})$,
$$
\sup_{t \in T} \left|\frac{1}{N} \sum_{i=1}^N \inr{X_i,t}^2 - \E \inr{X,t}^2 \right| \leq c_2\left(u^2d_2(T) \frac{E(T)}{\sqrt{N}} + u^4 \frac{E^2(T)}{N}\right),
$$
which improves the probability estimate from \cite{MenPao}.

\section*{Acknowledgements}
I am indebted to Vladimir Koltchinskii, Joe Neeman and Dong Xia for their careful reading of this manuscript and the many valuable comments and suggestions they have made.

\footnotesize {

\end{document}